\documentclass[a4paper]{article}%

\usepackage{amsmath,amsthm,amssymb,fancyhdr,color,epic,graphicx}%

\newtheorem{theorem}{\bf Theorem}[section]

\newtheorem{lemma}[theorem]{\bf Lemma}
\newtheorem{proposition}[theorem]{\bf Proposition}

\renewenvironment{proof}[1][Proof]{\noindent\textbf{#1.} }{\ \rule{0.4em}{0.4em}}
\newenvironment{acknowledgements}{\par\noindent\textbf{Acknowledgements.}}{\hfill\par}
\newenvironment{keywords}{\par\noindent\textbf{Keywords: }}{\hfill\par}

\theoremstyle{definition}
\newtheorem{algorithm}[theorem]{\bf Algorithm}
\newtheorem{remark}[theorem]{\bf Remark}

\definecolor{Be}{rgb} {0,0.08,0.45} 
\definecolor{darkgreen}{rgb}{0.00, 0.49, 0.00}
\definecolor{Mahog}{rgb}{0.6,  0.,   0}
\newcommand{\Rd}{\color{Mahog}}
\newcommand{\Be   }{\color{Be}}
\newcommand{\Bk}{\color{black}}
\newcommand{\Gn}{\color{darkgreen}}
\newcommand{\R}{I\!\!R}
\newcommand{\C}{\mathbb{C}}
\newcommand{\N}{I\!\!N}
\newcommand{\specialcell}[2][c]{\begin{tabular}[#1]{@{}c@{}}#2\end{tabular}}

\numberwithin{equation}{section}

\begin{document}

\title{A double scale fast algorithm for the transient evolution of a resonant tunneling diode}

\author{N. Ben Abdallah\footnote{Institut de Math\'ematiques de Toulouse (UMR CNRS 5219), Universit\'{e} Paul Sabatier, 31062
Toulouse Cedex 9, France}, A. Faraj\footnote{Institute of Natural Sciences, Shanghai Jiao Tong University, Shanghai 200240 P.R. China. Current address: Grenoble INP, ESISAR, 26902 Valence Cedex 9, France}}
\date{}

\maketitle

\begin{abstract}The simulation of the time-dependent evolution of the resonant tunneling diode is done by a multiscale algorithm exploiting the existence of resonant states. After revisiting the algorithm developed in [N. Ben Abdallah, O. Pinaud, J. Comput. Phys. 213, 1 (2006) 288-310] for the stationary case, the time-dependent problem is dealt with. The wave function is decomposed into a non resonant part and a resonant one. The projection method to compute the resonant part leads to an accurate algorithm thanks to a suitable interpolation of the non resonant one. The simulation times are largely reduced.
\end{abstract}

\section{Introduction}
The Resonant Tunneling Diode (RTD) is a heterostructure made of a superposition of thin layers of semiconductor of different types used to create double or multiple potential barriers. The possibility of resonant tunneling through multiple barriers and the resulting negative differential resistance were firstly discussed in \cite{Io,TsEs}. Such physical effects make the RTD useful in the conception of signal generators, detectors and mixers, multi-valued logic switches, low-power amplifiers, local oscillators, frequency locking circuits, and also as generators of multiple high frequency harmonics, see \cite{MiTa}.\\
In the mathematical and numerical models for the RTD, it is natural to incorporate the quantum effects described above. In particular, due to the small size of the device, the kinetic models for semi-conductors \cite{GrRa,Po} have to be replaced by their natural extension to the quantum level: the Wigner equation or equivalently a system of infinitely many Schr\"odinger equations. In the last models, the electron interaction can be represented thanks to a Poisson equation, which leads to the Wigner-Poisson \cite{BrMa,LiPa} or Schr\"odinger-Poisson system \cite{BeDeMa,BePi1,Ni}. Both approaches are able to recover the negative differential resistance in the I-V characteristics, however, it was shown in \cite{JLPrSj,PrSj} that the second approach could lead to hysteresis phenomena in agreement with physical observations.\\
For the Schr\"odinger-Poisson system, which will be adopted in the present work, the nonlinear coupling effects take place in the middle of the structure. This active region is out of equilibrium and connected to the exterior through access zones which allow the injection of electrons. Therefore, a first important issue is to find suitable transparent boundary conditions for the Schr\"odinger equation at the boundary of the active region. Such boundary conditions have been derived and analyzed in \cite{AnBe,BeDeMa,BePi1,Fr} and, at the discrete level, in \cite{Ar,EhAr,LeKi}.\\
Other theoretical studies were devoted to the Schr\"odinger-Poisson system, dealing with the well-posedness of the nonlinear problem \cite{BeDeMa,BePi1,Ni} and the derivation of simplified models for the RTD \cite{BoNiPa1,BoNiPa2,JLPrSj,PrSj}. Such simplified models were used, in the stationary case, in \cite{BoNiPa} to perform fast simulations, which agreement with the numerical solution of the complete Schr\"odinger-Poisson system was verified in \cite{BoFaNi}.\\
The simulation of the RTD is a challenging problem due to the presence of quantum resonances which create stiff variations of the wave function with respect to the frequency variable and increase dramatically the numerical complexity. In this context, the reduction of the simulation time provided by the simplified models described above is of particular interest. However, fast simulations with a higher accuracy were obtained in the stationary case in \cite{BePi2} by working on the complete Schr\"odinger-Poisson system. In addition to a new WKB finite element method (see also \cite{BeMoNe,Ne}), the reduction of the simulation time in this approach is provided by an adapted treatment of the transmission peaks using resonant modes.\\
In the time-dependent case, the supplementary numerical cost imposed by resonances is even more important and its reduction even more challenging. Indeed, the numerical complexity of the transparent boundary conditions of each time-dependent Schr\"odinger equation is of order one with respect to the number of time iterations, and the last is necessarily big since the period of the high energy free plane waves is much smaller than the scaling time of the system. Although advances (based on the adiabatic theorem) has been made in \cite{FaMaNi1,FaMaNi2} in order to find simplified models for the time-dependent simulations, an algorithm which extends the algorithm of \cite{BePi2} is still missing for the resolution of the transient Schr\"odinger-Poisson system in the presence of resonances. This is the aim of the present work.\\
In our RTD model, the classical observables, like electronic and current densities, are given by integrals on the frequency variable involving the wave function. To compute correctly these integrals, a uniform frequency mesh has to be very thin to be able to capture the resonance peaks. The high number of Schr\"odinger equations to be solved, and therefore the high numerical cost, of this standard approach was reduced in the stationary regime in \cite{Pi} using adaptive refinement. The last method does however not extend to the transient case since the resonant frequencies do move as time varies. Lately, inspired by an idea from \cite{PrSj}, the algorithm in \cite{BePi2} was proposed for the discretization of the stationary system. The method consists in decomposing the wave function in a resonant part living inside the double barrier and a non resonant part which is mostly localized outside. The non resonant part is smooth with respect to the frequency variable whereas the resonant one has sharp peaks at resonant frequencies. The latter is computed by a projection method after a precomputation of resonant states. Due to the sharpness of the resonance peaks, this process can not succeed without an accurate value of resonance energies and widths. Therefore, we first present an improvement of the method used in \cite{BePi2} to compute these values. The problem of finding quantum resonances has interested many authors, e.g. \cite{BiZw,RiMe}. In the one dimensional case, it can be written as a differential equation with homogeneous transparent boundary conditions selecting outgoing functions, see \cite{Fa,Pa}. Starting from this approach, we write the problem of finding resonances at the discrete level. Since the spectral parameter is involved nonlinearly in the boundary conditions, we obtain a nonlinear eigenvalue problem for which we derive a Newton-like method which extends to holomorphic matrices the method proposed in \cite{PeWi}. Additional references concerning the resolution of nonlinear eigenvalue problems can be found in \cite{Gui}.\\
In the present article, the Fermi-Dirac distribution involved in the densities operates like a frequency cutoff and only the first mode is important in the resonant part of the wave function. The multi-mode approximation follows the same line and finds applications in the simulation of the stationary two-dimensional RTD, see \cite{BeWu}.\\
The key point in the algorithm of \cite{BePi2} for the reduction of the simulation time is to use two scales for the frequency mesh. For the non-resonant part, the frequency mesh can be chosen coarse enough thus reducing the number of Schr\"odinger equations to be solved. For the resonant part, a single spatial problem (the computation of the resonant mode) is solved for all the frequency points. Therefore a refined frequency mesh can be used for the projection on the resonant mode, capturing this way the resonant peak without increasing significantly the numerical cost. Due to the scale difference between the two parts of the wave function, the projection method requires an interpolation of the non resonant part. In the present work, we show that the time-dependent algorithm is accurate only if the interpolation takes into account the time oscillations coming from the Schr\"odinger equation. In particular, we verify that the constant interpolation, which is used in the stationary case, is no more adapted in the time-dependent case.\\
For stability reasons, we solve the time-dependent Schr\"odinger equation with a semi-implicit scheme. In our nonlinear framework, such a scheme requires an adapted extrapolation of the potential at half-time step, see \cite{ChJiSu}, leading to an Adams-Bashforth-Crank-Nicolson finite difference method. We show that such an extrapolation allows bigger time steps than the extrapolation introduced in \cite{Pi} for our Schr\"odinger-Poisson system and also considered in \cite{Be} for a different kind of nonlinearity. In the stationary case, the nonlinearity of the Schr\"odinger-Poisson system is dealt with by using the self-consistent scheme in \cite{Gum}, see also \cite{Pi}. More numerical methods for the time-dependent nonlinear Schr\"odinger equation can be found in \cite{BaJiMa,JiWuYa}.\\
Together with the Adams-Bashforth-Crank-Nicolson scheme, we use the discrete boundary conditions in \cite{Ar,EhAr} for the time-dependent Schr\"odinger equation since they have the advantage to be verified exactly by the whole space discrete solution. We note that these boundary conditions can be used to write several kind of finite difference schemes, including splitting methods, to solve time-dependent nonlinear Schr\"odinger equations, see \cite{ZiEh}. Applying additional results from \cite{Ar,EhAr} to our RTD model, we verify numerically that, compared to the results in \cite{Pi}, the transparent boundary conditions for the Schr\"odinger equation can be fitted to improve the accuracy or simplified to reduce the computational time.\\
In order to be able to compare the time-dependent solution to the stationary one in a meaningful manner, the discretization of the space variable must be the same for the two regimes. Therefore, we will employ a standard finite difference method along with the stationary boundary conditions in \cite{Ar} to solve the stationary Schr\"odinger equation and concentrate on the frequency variable to handle the resonance peaks. Similarly, in our algorithm, the resonant part is computed using finite difference in order to be consistent with the non resonant one. For that purpose, we write the discrete boundary conditions for the resonant mode as an extension to complex valued energies of the homogeneous stationary boundary conditions. These boundary conditions are verified exactly by the resonant mode thanks to the following strategy, which is similar to the strategy in \cite{Ar,EhAr}: first we consider the discretization of the differential equation on the whole space, and then we derive the boundary conditions for the difference scheme directly on a purely discrete level.\\
The paper is organized as follows. In section \ref{sec_model}, are presented the stationary and time-dependent Schr\"odinger-Poisson systems. The stationary case is dealt with in section \ref{sec_stat} where are given the reference method, the adaptive refinement method and two methods based on the projection on the resonant mode. Our method to compute resonances and resonant modes is also introduced. In section \ref{sec_instat}, the time-dependent reference method is recalled and our improved transient algorithm is presented. The numerical results are given in section \ref{sec_results}. The different methods are compared in the stationary and in the time-dependent case. In the time-dependent case, the comparison shows the improvement provided by our new algorithm and the existence of two resonant peaks is verified. In the appendix are given the finite difference methods with discrete transparent boundary conditions to solve the stationary and time-dependent Schr\"odinger equation. For the time-dependent scheme, a stability result is proven in the situation of a continuous influx in the device. The discrete problem for the computation of resonances is also obtained and the properties of its solution are presented.

\section{The model}\label{sec_model}

\begin{figure}[!ht]
\begin{center}\begin{picture}(230,90)
\put(5,30){\line(1,0){35}}
\put(40,30){\line(4,-1){20}}
\put(60,25){\line(0,1){70}}
\put(60,95){\line(4,-1){10}}
\put(70,92.5){\line(0,-1){70}}
\put(70,22.5){\line(4,-1){10}}
\put(80,20){\line(0,1){70}}
\put(80,90){\line(4,-1){10}}
\put(90,87.5){\line(0,-1){70}}
\put(90,17.5){\line(4,-1){20}}
\put(110,12.5){\line(1,0){40}}
\Be\put(55,25){\vector(0,1){70}} \put(40,55){$v_0$}
\put(115,30){\vector(0,-1){17.5}} \put(120,21){$-B(t)$}\Bk
\Gn \dottedline(40,30)(150,30) \Bk
\Rd \put(5,20){$0$} \put(150,2.5){$L$} \Bk
\put(35,22){$a_1$}\put(110,2.5){$b_1$}\put(52,18){$a_2$}\put(88,8){$b_2$}\put(64.5,16){$a_3$}\put(76.5,13){$b_3$}
\end{picture}\end{center}
\caption{External potential $U(t)$.}\label{figVe}
\end{figure}
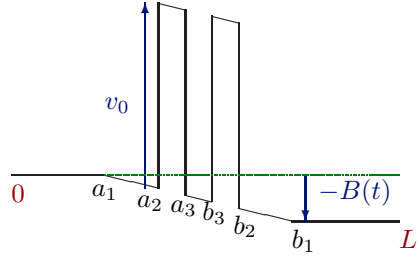

It consists in an infinite number of Schr\"odinger equations coupled to the Poisson equation. The Schr\"odinger equation involves the time-dependent Hamiltonian
\[
H(t)=-\frac{\hbar^2}{2m}\partial_x^2 + U(t) + V(t),
\]
where $\hbar$ is the reduced Planck constant, $m$ is the effective mass of the electron and $x\in \R$ is the position variable. The domain occupied by the device is the interval $[0,L]$, $L>0$. The self-consistent potential $V(t)$ is due to Coulomb interaction and depends nonlinearly on the wave functions. The external potential $U(t)$ describes the double barrier and the applied bias, see Figure \ref{figVe}, and is given as if a data of the problem. Given the points
\[
0 < a_1 < a_2 < a_3 < b_3 < b_2 < b_1 < L\,,
\]
the external potential writes
\[
U(t)= v_0{\bf1}_{[a_2,a_3[\cup]b_3,b_2]}-B(t)(\frac{x-a_1}{b_1-a_1}{\bf1}_{[a_1,b_1[}+{\bf1}_{[b_1,+\infty[})\,,
\]
where $v_0 \geq 0$ and $B(t) \geq 0$ are scalars representing respectively the height of the barrier and the amplitude of the applied bias. The points $a_1$ and $b_1$ are the extremities of the diode. The interval $[a_3,b_3]$ is the quantum well.
The problem is the following infinite system of Schr\"odinger equations
\begin{equation}\label{schInst}
\displaystyle \left\{
\begin{array}{ll}i\hbar\partial_t\Psi_k(t) = -\frac{\hbar^2}{2m}\partial_x^2\Psi_k(t) + (U(t)+V(t))\Psi_k(t)\,, & x\in\R\\[1.65mm]
\Psi_k(0) = \Phi_k
\end{array}\right.\,,
\end{equation}
where the initial condition $\Phi_k$ is a generalized eigenfunction of the initial Hamiltonian $H(0)$ corresponding to the frequency $k\in\R$. The functions $\Phi_k$ are defined more in details at the end of this section. The self-consistent potential satisfies the Poisson equation:
\begin{equation}\label{eqPot}
\displaystyle\left\{
\begin{array}{ll} -\partial_x^2V(t) = \frac{q^2}{\varepsilon}(n(t)-n_D)\,, & x\in(0,L)\\[1.65mm]
V(t,0)=V(t,L)=0
\end{array}\right.
\end{equation}
and is extended by $0$ outside $(0,L)$. In \eqref{eqPot}, $q$ is the elementary charge of the electron, $\varepsilon$ is the dielectric constant, $n_D$ is the doping equal to 
\[
\displaystyle n_D = n_D^1({\bf 1}_{[0,a_1[} + {\bf 1}_{]b_1 ,L]}) + n_D^2{\bf 1}_{[a_1,b_1]}\,,
\] 
with $n_D^1 > n_D^2 \geq 0$. The electron density $n(t)$ is given by:
\begin{equation}\label{densIns}
\displaystyle n(t,x) = \int_{\R}g(k)|\Psi_k(t,x)|^2dk\,.
\end{equation}
In our model, the injection profile $g$ is the one dimensional Fermi-Dirac distribution:
\begin{equation}\label{FermiDirac}
\displaystyle g(k) = \frac{mk_BT}{2\pi^2\hbar^2}\ln\left(1 + \exp\left(\frac{E_F-\frac{\hbar^2k^2}{2m}}{k_BT}\right)\right)\,,
\end{equation}
where $k_B$ is the Boltzmann constant, $T$ is the temperature of the semiconductor and $E_F$ is the Fermi level.\\ 
In the stationary regime, the external potential is the potential $U_I=U(0)$ corresponding to the initial bias $B_I=B(0)$. The time-independent problem is the following infinite system of equations
\begin{equation}\label{VpGe}
\displaystyle -\frac{\hbar^2}{2m}\partial_x^2\Phi_k + (U_I+V_I)\Phi_k = E_k\Phi_k\,, \quad x\in\R\,,
\end{equation}
with scattering conditions
\begin{equation}\label{eqSCkp}
\displaystyle \left\{\begin{array}{ll} \Phi_k(x) = e^{ikx} + r(k) e^{-ikx}\,,& x<0\\[1.65mm]
\Phi_k(x) = t(k)e^{i\sqrt{k^2+2mB_I/\hbar^2}x}\,, &x > L
\end{array}\right.
 \quad  \textrm{for } k \geq 0
\end{equation}
and
\begin{equation}\label{eqSCkn}
\displaystyle \left\{\begin{array}{ll} \Phi_k(x) = t(k) e^{-i\sqrt{k^2-2mB_I/\hbar^2}x}\,,& x<0\\[1.65mm]
\Phi_k(x) = e^{ikx} + r(k) e^{-ikx}\,, &x > L
\end{array}\right. \quad  \textrm{for } k < 0\,,
\end{equation}
where 
\begin{equation}\label{eqDisp}
\displaystyle E_k=\left\{\begin{array}{ll} \frac{\hbar^2k^2}{2m}\,, & k\geq 0\\[1.65mm]
\frac{\hbar^2k^2}{2m}-B_I\,, & k<0\end{array}\right.\,.
\end{equation}
 As in the time-dependent case, the self-consistent potential $V_I$ satisfies the Poisson equation:
\begin{equation}\label{eqPotst}
\displaystyle\left\{
\begin{array}{ll} -\partial_x^2V_I = \frac{q^2}{\varepsilon}(n_I-n_D)\,, & x\in(0,L)\\[1.65mm]
V_I(0)=V_I(L)=0
\end{array}\right.
\end{equation}
and is extended by $0$ outside $(0,L)$. The electron density $n_I$ is given by:
\begin{equation}\label{denStat}
\displaystyle n_I(x) = \int_{\R}g(k)|\Phi_k(x)|^2dk\,.
\end{equation}

\section{The stationary algorithm revisited }\label{sec_stat}

\subsection{Recalling the standard algorithm}\label{sec_algStat}

To solve the stationary nonlinear Schr\"odinger-Poisson system \eqref{VpGe}\eqref{eqPotst}\eqref{denStat}, we will use the same algorithm as in \cite{Pi}. It is based on a Gummel iteration, see \cite{Gum}, which consists in the computation of a sequence $\left(V_I^l\right)_{l\in\N}$ where the potential $V_I^{l+1}$ at step $l+1$ is deduced form the potential $V_I^l$ at step $l$ by solving the following nonlinear equation:
\begin{equation}\label{itGum}
\displaystyle \left\{
\begin{array}{ll}-\partial_x^2V_I^{l+1} = \frac{q^2}{\varepsilon}(n_I^l\exp((V_I^l-V_I^{l+1})/V_{ref})-n_D)\,, &x\in(0,L)\\[1.65mm]
V_I^{l+1}(0)=V_I^{l+1}(L)=0
\end{array}\right.
\end{equation}
for a fixed reference potential $V_{ref}$. In \eqref{itGum}, $n_I^l$ is the density \eqref{denStat} where the wave functions are solution to \eqref{VpGe} with $V_I=V_I^l$.\\
The repartition function $g$ being exponentially decreasing at infinity, the integral in \eqref{denStat} can be restricted, in computations, to a domain $[-\kappa,\kappa]$ where $\kappa$ is chosen to be big enough. We consider a discretization $x_0=0,x_1,...,x_j,...,x_J=L$ of the interval $[0,L]$ with uniform grid spacing $\Delta x$ and a discretization $k_0=-\kappa,k_1,...,k_p,...,k_P=\kappa$ of the interval $[-\kappa,\kappa]$ with non-necessarily uniform grid spacing $\Delta k_{p}=k_{p+1}-k_p$. In order to verify, the condition \eqref{eqstabx} for all $k\in[-\kappa,\kappa]$, the grid spacing for the space variable is chosen such that:
\begin{equation}\label{eqstabx_toutk}
\Delta x < \frac{1}{\sqrt{\kappa^2+\frac{2m}{\hbar^2}B_I}}\,.
\end{equation}
Then, for a given initial potential $V_I^0$, the algorithm writes:
\begin{algorithm}[Gummel algorithm]\label{algGum}
\hfill\par\noindent
Fix $l=0$.\\
Set $\Delta = \textrm{tol} + 1$.\\
\underline{Do While} $\Delta\geq\textrm{tol}$:
\begin{itemize}  
\item[]{\bf Computation of the density}
\begin{itemize} 
\item[$S1.$]For $p=0,...,P$: computation of the wave function $\Phi^l_p$ at the frequency $k_p$ from the potential $V_I^l$ by solving equation \eqref{VpGe} with $V_I=V_I^l$.
\item[$S2.$]Numerical integration: for $j=0,...,J$
\begin{equation}\label{eqtrap}
\displaystyle n^l_{I,j}= \sum_{p=0}^{P-1}\left(g(k_p)\left|\Phi^l_{p,j}\right|^2+g(k_{p+1})\left|\Phi^l_{p+1,j}\right|^2\right)\frac{\Delta k_{p}}{2}\,.
\end{equation}
\end{itemize} 
\item[]{\bf Coupling to the Poisson equation :  Gummel iteration}\\
Computation of the potential $V_I^{l+1}$ from the potential $V_I^l$ and the density $n_I^l$ by solving equation \eqref{itGum}.
\end{itemize} 
Set 
\begin{equation}\label{eq_erreur}
\Delta = \frac{\left\vert V_I^{l+1}-V_I^l\right\vert_2}{\left\vert V_I^{l+1}\right\vert_2}\,.
\end{equation}
Set $l=l+1$.\\
\underline{End Do}.
\end{algorithm}
In equation \eqref{eqtrap}, $n^l_{I,j}$ denotes the approximation of the density $n^l_I(x_j)$ and $\Phi^l_{p,j}$ the approximation of the wave function $\Phi^l_p(x_j)$. In equation \eqref{eq_erreur} and in what follows, $\vert v \vert_2$ denotes the $l^2$-norm of the vector $v\in\C^{J+1}$. Equation (\ref{itGum}) is nonlinear and is solved with a Newton method where the Laplacian is discretized using finite difference at the points $(x_j)_{0\leq j \leq J}$.\\
Because of the peaked form of the transmission near resonances, the method used to compute the density $n_I^l$ from the potential $V_I^l$ (steps $S1$ and $S2$) is of major importance, and the Gummel algorithm may fail to converge if a non accurate method is used.\\
In the case of the finite difference discretization, the strategy of \cite{Pi} can be described as follows. The step $S1$ is performed without particular treatment of resonances: the function $\Phi^{l}_p$ is computed on $[0,L]$ by the resolution of the problem \eqref{VpGe}-\eqref{eqDisp} with the finite difference scheme \eqref{eqdfk} and the discrete transparent boundary conditions \eqref{eqBCkpg}, \eqref{eqBCkpd} for $k\geq 0$ and \eqref{eqBCkng}, \eqref{eqBCknd} for $k<0$ where the potential $Q_I$ is replaced by $U_I+V^l_I$. In that case, the convergence is provided by the choice of the frequency mesh $\{k_p,\,p=0,...,P\}$ for the trapezoidal rule in the step $S2$. The reference method is obtained by taking the grid spacing $\Delta k_{p}$ equal to a constant $\Delta k>0$ small enough. This method will be called Direct Resolution. In order to reduce the high numerical cost, which is due to the big number of Schr\"odinger equations to be solved, an Adaptive Method is proposed in \cite{Pi}. For the Adaptive Method, the grid spacing $\Delta k_{p}$ is taken equal to the small constant $\Delta k$ only when $E_{k_p}$ is close to a resonant energy, otherwise $\Delta k_{p}=\nu\Delta k$ where $\nu\geq 2$. The frequency points where the mesh is refined are detected using the logarithmic derivative, with respect to the frequency $k$, of the transmission coefficient related to the wave function $\Phi_k$ (see \cite{Pi,Fa}).\\
As noticed in the introduction, the Adaptive Method can not be generalized to the time-dependent case. In this context, the One Mode Approximation algorithm presented in section \ref{sec_1mSt} is very useful. As proposed in \cite{BePi2}, this algorithm consists in decomposing the wave function in a non resonant part and a resonant one proportional to the first resonant mode. In \cite{BePi2}, using a WKB interpolation to compute each part of the wave function, an adapted treatment of the step $S1$ is realized.\\
In the present work, we do not use the WKB interpolation, however the accuracy required at step $S1$ is reached by an improvement of the computation of the first resonance. In particular, a precise computation of its imaginary part is essential. A simple reconstitution of the wave functions not allowing to modify the frequency mesh, the gain in computational time is reached by adapting the step $S2$. In the present work, two strategies are proposed for the step $S2$. The first method is inspired from \cite{BePi2} and consists in using different frequency scales for the different parts of the wave function. The frequency mesh is taken thin only for the coefficient of proportionality to the resonant mode and the number of Schr\"odinger equations to be solved is reduced. This method will be denoted One Mode Approximation and will be generalized to the time-dependent case. A variant of this method is obtained by an explicit integration of the coefficient of proportionality to the resonant mode instead of the trapezoidal rule. In that case all the resonant peak is integrated and the frequency mesh can be chosen coarse enough. This method will be denoted One Mode Approximation Integ.\\
The algorithms Direct Resolution, Adaptive Method, One Mode Approximation and One Mode Approximation Integ are compared in section \ref{secComp}. 

\subsection{Accurate computation of resonances}\label{sec_calres}

In the founder work of \cite{AgCo}, see also \cite{HiSi}, resonances of a self-adjoint operator are defined using analytic transformations, and it is a common fact that it corresponds to an eigenvalue in a modified $L^2$ space \cite{FuRa}. Using the second approach, the authors in \cite{Pa} and \cite{Fa} define resonances for Hamiltonians of the form
\begin{equation}\label{eqHR}
\displaystyle -\frac{\hbar^2}{2m}\partial_x^2+Q\,,
\end{equation}
where the potential $Q$ verifies
\begin{equation}\label{eqHQN}
\displaystyle Q(x)=0\,,\, x\leq 0 \quad \textrm{ and } \quad  Q(x)=Q_L\,, \, x\geq L \quad \textrm{ where } \quad Q_L\leq 0\,.
\end{equation}
Starting from this definition, the resonances of the Hamiltonian \eqref{eqHR} can be characterized as the numbers $z$ in $\mathbb{C}\setminus\{i\R\cup i\R+Q_L\}$ such that the problem
\begin{equation}\label{pbMres}
\displaystyle \left\{\begin{array}{lr}[-\frac{\hbar^2}{2m}\partial_x^2+Q]u = zu\,,& x\in\R\\[1.65mm]
||u||_{L^2(0,L)}=1
\end{array}\right.
\end{equation}
has solutions which are purely outgoing outside the interval $[0,L]$. The function $u$ in \eqref{pbMres} is called resonant mode associated to the resonance $z$. As described in Appendix \ref{app_RDF}, the eigenproblem \eqref{pbMres} is restricted (for numerical computations) to the domain $[0,L]$ by the aid of transparent boundary conditions selecting only outgoing states. In particular, this problem is discretized using the finite difference scheme \eqref{eqdfr} inside the domain and the boundary conditions \eqref{eqBCgr}, \eqref{eqBCdr} at $x=0$ and $x=L$. Since the boundary conditions depend nonlinearly on the spectral parameter $z$, this process leads to the following nonlinear eigenvalue problem: find $(u,z)\in\C^{J+1}\times\C$ such that
\begin{equation}\label{pbMresDisc}
\displaystyle \left\{\begin{array}{l}M(z)u=0\\[1.65mm]
u^Hu=1
\end{array}\right.\,,
\end{equation}
where $M(z)$ is the nonlinear matrix valued function of $z$ given by \eqref{eqMz} and $u^H$ denotes the complex conjugate transpose of the vector $u$. A Newton-like method is introduced in \cite{PeWi} to solve \eqref{pbMresDisc} when $M(z)$ is a polynomial. Following the same line, we derive a method for the general case where $M(z)$ is holomorphic. Note first that the term $u^Hu$ appearing in \eqref{pbMresDisc} is not differentiable with respect to $u$, therefore a Newton method must be modified to be applied here. For a given iterate $(u^n,z^n)$ verifying $(u^n)^Hu^n=1$, we are looking for a direction $(\delta u^n,\delta z^n)$ such that $(u^n+\delta u^n,z^n+\delta z^n)$ is solution to the problem \eqref{pbMresDisc}. Using $(u^n)^Hu^n=1$ the system:
\begin{equation*}
\displaystyle \left\{\begin{array}{l}M(z^n+\delta z^n)(u^n+\delta u^n)=0\\[1.65mm]
(u^n+\delta u^n)^H(u^n+\delta u^n)=1
\end{array}\right.
\end{equation*}
gives at order 2:
\begin{equation}\label{eqNorm}
\displaystyle 
\left\{\begin{array}{l}M(z^n)\delta u^n+\delta z^nM'(z^n)u^n = -M (z^n)u^n\\[1.65mm] 
(u^n)^H \delta u^n + (\delta u^n)^Hu^n = 0
\end{array}\right.\,.
\end{equation}
As noticed earlier, the second equation in \eqref{eqNorm} does not have the suitable form. To deal with this problem, we remark that it is enough to impose $(u^n)^H \delta u^n=0$ in order to verify this equation. This leads to the following linear system:
\begin{equation}\label{newtRes}
\displaystyle \left[\begin{array}{cc}M(z^n) & M'(z^n)u^n \\ (u^n)^H & 0\end{array}\right]
\left[\begin{array}{c}\delta u^n \\ \delta z^n\end{array}\right]
= \left[\begin{array}{c}-r^n\\0\end{array}\right]\,, \quad \textrm{ where }  r^n=M(z^n)u^n\,.
\end{equation}
Its resolution corresponds to an iteration of our Newton-like method to compute the resonant mode and the resonance. The assumption $(u^n)^Hu^n=1$, that we made to obtain the system \eqref{newtRes}, is verified at order $2$ as long as $(u^0)^Hu^0=1$.

\subsection{The One Mode Approximation}\label{sec_1mSt}

We start with the description of the step $S1$ of the computation of a wave function $\Phi_k$ solution to the stationary Schr\"odinger equation \eqref{VpGe} for a given frequency $k$. To simplify the notations, the exponent $l$ appearing in Algorithm \ref{algGum} will be omitted. Following the works \cite{BePi2} and \cite{PrSj}, the One Mode Approximation consists in the decomposition:
\begin{equation}\label{form1mStat1}
\displaystyle \Phi_k=\Phi^{nr}_k+\Phi^r_k\,,
\end{equation}
where the non resonant part $\Phi^{nr}_k$ solves the stationary Schr\"odinger equation:
\begin{equation}\label{SchroStatNR}
\displaystyle [-\frac{\hbar^2}{2m}\partial_x^2+U_{I,fill}+V_I]\Phi^{nr}_k = E_k\Phi^{nr}_k\,, \quad x\in\R\,,
\end{equation}
with filled potential $U_{I,fill}=U_I+v_0{\bf1}_{[a_3,b_3]}$. Equation \eqref{SchroStatNR} comes with the scattering conditions \eqref{eqSCkp}, \eqref{eqSCkn} and the relation \eqref{eqDisp}. It is solved on $[0,L]$ using the finite difference scheme \eqref{eqdfk} and the discrete transparent boundary conditions \eqref{eqBCkpg}, \eqref{eqBCkpd} for $k\geq 0$ and \eqref{eqBCkng}, \eqref{eqBCknd} for $k<0$ where the potential $Q_I$ is replaced by $U_{I,fill}+V_I$. The statistic $g$ being rapidly decreasing at infinity, the resonant part $\Phi^r_k$ is sought, on $[0,L]$, proportional to the resonant mode $u_I$ of minimal resonant energy $\textrm{Re}\,(z_I)$ solution to 
\begin{equation}\label{eqResI}
[-\frac{\hbar^2}{2m}\partial_x^2+U_I+V_I]u_I=z_Iu_I \quad \textrm{ and } \quad \int_0^L|u_I(x)|^2dx=1\,.
\end{equation}
In other words, we look for $\Phi^r_k$ of the form:
\begin{equation}\label{form1mStat2}
\displaystyle \Phi^r_k(x) =\theta_k u_I(x), \quad x\in[0,L] \,.
\end{equation}
Inserting \eqref{form1mStat1} and \eqref{form1mStat2} in the stationary Schr\"odinger equation \eqref{VpGe}, and using \eqref{SchroStatNR} and \eqref{eqResI}, we get the following explicit value of the proportionality coefficient:
\begin{equation}\label{eqTheta}
\displaystyle \theta_k = \frac{1}{z_I-E_k}v_0\int_{a_3}^{b_3}\Phi^{nr}_k\overline{u_I}dx\,.
\end{equation}
The resonant mode and the resonance are computed using the method presented in section \ref{sec_calres} with the potential $Q=U_I+V_I$. This method is initialized at the fundamental mode and fundamental energy of the Hamiltonian: 
\begin{equation}\label{hamBar}
\displaystyle [-\frac{\hbar^2}{2m}\partial_x^2+U_I+V_I]
\end{equation}
equipped with homogeneous Dirichlet boundary conditions at $a_2$ and $b_2$. It is shown in \cite{BoNiPa1,BoNiPa2} that the real part of the resonances are close to the eigenvalues of the Dirichlet Hamiltonian \eqref{hamBar}. The imaginary part of resonances being small, such an initialization insures the convergence of the algorithm to the resonance with smaller energy. This achieves the step $S1$.\\
For the One Mode Approximation, the step $S2$ consists in using the wave function decomposition \eqref{form1mStat1} and \eqref{form1mStat2} in order to compute accurately the density with a low numerical cost. First, using an argument of localization of support, see \cite{PrSj}, the cross term in the development of $\left|\left(\Phi^{nr}_k+\theta_ku_I\right)(x)\right|^2$ can be neglected, which gives the approximation:
\begin{equation}\label{eq_densresst}
\displaystyle n_I(x)=\int_{\R}g(k)\left|\Phi^{nr}_k(x)\right|^2dk+\int_{\R}g(k)\left|\theta_k\right|^2dk\left|u_I(x)\right|^2\,.
\end{equation}
The non resonant part of the wave function being regular with respect to $k$, the first integral above can be computed using a trapezoidal rule with coarse frequency mesh, reducing this way the number of Schr\"odinger equations to be solved. As remarked above, the resonance $z_I$ has a small imaginary part. Therefore, it appears from formula \eqref{eqTheta}, that the coefficient $\theta_k$ has a sharp peak when the wave function energy $E_k$ is close to the resonant energy $\textrm{Re}\,(z_I)$. The key-point in the step $S2$ is then the approximation of the integral $\int_{\R}g(k)\left|\theta_k\right|^2dk$.\\
For the sake of simplicity, the non resonant part of the wave function will be called non resonant wave function as well. The same remark holds for the time-dependent One Mode Approximation presented in section \ref{sec_1mInst}.  

\subsubsection{First method: One Mode Approximation}\label{sec1mapp}

Consider a refined discretization $k_0=-\kappa,k_1,...,k_p,...,k_P=\kappa$ of the interval $[-\kappa,\kappa]$ with uniform grid spacing $\Delta k$ and suppose that there exist two integers $P'$ and $\nu$ such that $\frac{P}{P'}=\nu$. Consider the large scale mesh $\hat{k}_{p'}=k_{\nu p'}$ for $p'=0,...,P'$ with uniform grid spacing $\Delta \hat{k}=\nu\Delta k$. Then, if we note $\Phi^{nr}_{p',j}$ the approximation of the non resonant wave function $\Phi^{nr}_{k_{\nu p'}}(x_j)$, $u_{I,j}$ the approximation of the resonant mode $u_I(x_j)$ and $\theta_p$ the approximation of the coefficient $\theta_{k_p}$, the formula for the density is
\begin{multline*}
n_{I,j}=\sum_{p'=0}^{P'-1}\left(g(\hat{k}_{p'})\left|\Phi^{nr}_{p',j}\right|^2+g(\hat{k}_{p'+1})\left|\Phi^{nr}_{p'+1,j}\right|^2\right)\frac{\Delta \hat{k}}{2}\\ +\sum_{p=0}^{P-1}\left(g(k_{p})\left|\theta_p\right|^2+g(k_{p+1})\left|\theta_{p+1}\right|^2\right)\frac{\Delta k}{2}\left|u_{I,j}\right|^2\,.
\end{multline*}
The crucial point for the reduction of the simulation time is that only $P'$ Schr\"odinger equations have to be solved, instead of $P$ equations for the Direct Resolution. However, this reduction implies that we only have access to the functions $\Phi^{nr}_{p',j}$ for $p'=0,...,P'$ and the computation of the coefficient $\theta_{p}$, for $p\neq \nu p'$, requires an interpolation of the non resonant wave function in the r.h.s. of equation \eqref{eqTheta}. It is provided by the piecewise constant interpolation formula below: for $p=0,...,P$
\begin{equation}\label{eqThetaIp}
\displaystyle \theta_{p} = \frac{1}{z_I-E_{p}}v_0\sum_{j\in w}\Phi^{nr}_{p',j}\overline{u_{I,j}}\Delta x\,,
\end{equation}
where $w=\{j\,|\,a_3\leq x_j<b_3\}$, $E_{p}=E_{k_p}$ and $0\leq p'\leq P'$ is the integer such that $p=\nu p'+r$ for some $0\leq r\leq \nu-1$. Since the evaluation of the coefficient $\theta_p$ in formula \eqref{eqThetaIp} is numerically cheap, the number $P$ can be chosen big enough to catch the resonance peak without a significant increase of the computational cost.

\subsubsection{Second method: One Mode Approximation Integ}\label{sec1mInt}

We only considered here the large scale frequency mesh $\hat{k}_{0},...,\hat{k}_{P'}$ introduced in section \ref{sec1mapp}. As noticed above, the integral of the non resonant part can be computed with a simple trapezoidal rule. Then, if we set $J_p:=\int_{\hat{k}_{p}}^{\hat{k}_{p+1}}g(k)\left|\theta_k\right|^2dk$, we get the following formula for the density:
\begin{equation}\label{eqAppnj}
\displaystyle n_{I,j}=\sum_{p=0}^{P'-1}\left(g(\hat{k}_{p})\left|\Phi^{nr}_{p,j}\right|^2+g(\hat{k}_{p+1})\left|\Phi^{nr}_{p+1,j}\right|^2\right)\frac{\Delta \hat{k}}{2}
+\sum_{p=0}^{P'-1}J_p\left|u_{I,j}\right|^2\,.
\end{equation}
In that case, the required precision for the density is reached thanks to an adapted treatment of the integral in $J_p$. It follows from \eqref{eqTheta} that:
\[
\displaystyle J_p=\int_{\hat{k}_{p}}^{\hat{k}_{p+1}}S_k\frac{k}{\left|E_k-z_I\right|^2}dk\,,
\]
where
\[
\displaystyle S_k=\frac{g(k)}{k}v_0^2\left|\int_{a_3}^{b_3}\Phi^{nr}_k\overline{u_I}dx\right|^2\,.
\]
Using numerical observations, it is shown in \cite{Fa} that $\Vert\Phi^{nr}_k\Vert_{L^{\infty}(0,L)}=\mathcal{O}(k)$ when $k$ tends to $0$. This limit is also verified theoretically in the same work, under the condition that $k>0$ and for particular potentials. Therefore, $S_k$ is a smooth function of the variable $k$ in a small neighborhood of $0$. Since the wave function $\Phi^{nr}_k$ has no resonance peak, $S_k$ is also smooth outside such a neighborhood. As a consequence, $S_k$ can be approximated by a constant on $[\hat{k}_{p},\hat{k}_{p+1}]$ which gives:
\[
\displaystyle J_p=S_p\int_{\hat{k}_{p}}^{\hat{k}_{p+1}}\frac{k}{\left|E_k-z_I\right|^2}dk\,,
\]
where 
\[
S_p:=\frac{g\left(\hat{k}_{p}\right)}{\hat{k}_{p}}v_0^2\left|\sum_{j\in w}\Phi^{nr}_{p,j}\overline{u_{I,j}}\Delta x\right|^2
\]
and $w$ is defined in section \ref{sec1mapp}. The resonance $z_I$ can be written 
\begin{equation}\label{eqZI}
z_I=E_I-i\Gamma_I\,,
\end{equation}
where $E_I\in\R$ and $\Gamma_I>0$. Under the assumption $0 \notin [\hat{k}_{p},\hat{k}_{p+1}]$, it holds $E_k=\frac{\hbar^2k^2}{2m}-B_p$ on $[\hat{k}_{p},\hat{k}_{p+1}]$, where $B_p$ is a constant equal to $0$ or $B_I$. Then, it follows:
\[
\displaystyle J_p=S_p\int_{\hat{k}_{p}}^{\hat{k}_{p+1}}\frac{k}{\left(\frac{\hbar^2k^2}{2m}-B_p-E_I\right)^2+\Gamma_I^2}dk\,.
\]
Utilizing the change of variable $E=\frac{\hbar^2k^2}{2m}-B_p-E_I$ in the above integral, we obtain the following approximation of $J_p$:
\begin{equation}\label{eqAppJp}
\displaystyle J_p=\frac{m S_p}{\hbar^2\Gamma_I}\left(\arctan\left(\frac{\hat{E}_{p+1}-E_I}{\Gamma_I}\right)-\arctan\left(\frac{\hat{E}_{p}-E_I}{\Gamma_I}\right)\right) \quad \textrm{ for  } \quad p=0,...,P'\,,
\end{equation}
where $\hat{E}_p:=E_{\hat{k}_{p}}$. We will use \eqref{eqAppJp} regardless of the verification of the assumption $0 \notin [\hat{k}_{p},\hat{k}_{p+1}]$. It provides an accurate treatment of the resonance peak even if $\Delta \hat{k}$ is large. The approximation of the density follows by replacing \eqref{eqAppJp} in \eqref{eqAppnj}.
\begin{remark}\label{rqInitres}
The computational time of the One Mode Approximation algorithm can be reduced when the method to compute the resonance is initialized at the fundamental mode and fundamental energy  $(u_D,E_D)$, of the Dirichlet Hamiltonian \eqref{hamBar}, only for the first iteration $l=0$ of the Gummel algorithm \ref{algGum}. Indeed, for $l\geq 1$, the computation of $(u_D,E_D)$ is avoided by initializing the method at the resonant mode and resonance obtained at the previous iteration $l-1$.  
\end{remark}
\begin{remark}
The size of the grid spacing $\Delta k$ for the refined frequency mesh is chosen numerically. It can be estimated analytically by imposing that, around the resonance energy $E_I$, the difference $E_{p+1}-E_{p}$ is smaller than half the resonance peak width $\Gamma_I$.  
\end{remark}

\section{The time-dependent algorithm}\label{sec_instat}

This section deals with the time evolution of the system under a time-dependent applied bias $B(t)$ as in Figure \ref{figVe}. We will consider only the case of the step like bias of the form:
\begin{equation}\label{B_marche}
\displaystyle B(t)=\left\{\begin{array}{ll} B_I\,, & t= 0\\[1.5mm]
B_{\infty}\,, & t>0\end{array}\right.\,,
\end{equation}
where $B_I, B_{\infty} \geq 0$ and $B_I\neq B_{\infty}$. The case of the general bias requires an adaptation of the discrete transparent boundary conditions for the resolution of the transient Schr\"odinger equation.

\subsection{The algorithm}\label{sec_algInst} 

As in the stationary regime, we first present the algorithm proposed in \cite{Pi} to solve the time-dependent nonlinear Schr\"odinger-Poisson system \eqref{schInst}\eqref{eqPot}\eqref{densIns}. Then, we give the details of our time-dependent One Mode Approximation algorithm.\\
The space and frequency meshes are defined as in section \ref{sec_algStat}. In order to verify the condition \eqref{eqstabt} for all $k\in[-\kappa,\kappa]$, the time step $\Delta t$ is chosen such that:
\begin{equation}\label{eqstabt_toutk}
\Delta t < \frac{\pi\hbar}{2\left(\frac{\hbar^2\kappa^2}{2m}+\sup_{t\geq0}B(t)\right)}\,.
\end{equation}
Then, for a given exterior potential $U(t)$, the algorithm corresponds to the computation of the sequence $V^l$ of approximations of the self-consistent potential at time $t^l=l\Delta t$. It provides a sequence of densities $n^l$ at time $t^l$. The initial potential $V^0$, resp. wave function $\Psi^0_p$, is given by $V_I$, and resp. $\Phi_{k_p}$, solution to the stationary Schr\"odinger-Poisson system \eqref{VpGe}\eqref{eqPotst}\eqref{denStat} with exterior potential equal to $U_I$. The Schr\"odinger equation \eqref{schInst} is solved using a Crank-Nicolson scheme and, in our nonlinear framework, a choice of the intermediary  potential $V^{l+\frac{1}{2}}$ which preserves the unconditional stability of the semi-implicit scheme is required. In particular, it will be shown in section \ref{sec_parInst} that the extrapolation:
\begin{equation}\label{eqExtV}
V^{l+\frac{1}{2}}=\frac{3}{2}V^l-\frac{1}{2}V^{l-1}
\end{equation}
allows bigger time steps than the following extrapolation proposed in \cite{Pi}:
\begin{equation}\label{eqExtV_pi}
V^{l+\frac{1}{2}}=2V^l-V^{l-\frac{1}{2}}\,.
\end{equation}
The extrapolation \eqref{eqExtV} is proposed in \cite{ChJiSu} and corresponds to a modified Adams-Bashforth-Crank-Nicolson method. We can now give the time-dependent algorithm:
\begin{algorithm}[Transient algorithm]\label{AlgTra}
\hfill\par
\underline{Do For} \, $l=0$ \, \underline{to} \, $l_{tot}-1$:
\begin{itemize}  
\item[]{\bf Computation of the intermediary  potential}\\
Computation of the potential $V^{l+\frac{1}{2}}$ from the potentials $V^l$ and $V^{l-1}$ using \eqref{eqExtV} with the convention $V^{-1}=V^0$. 
\item[]{\bf Computation of the density}
\begin{itemize} 
\item[$S1.$]For $p=0,...,P$: computation of the wave function $\Psi^{l+1}_p$, solution to \eqref{schInst} at time $t^{l+1}$ and frequency $k_p$, using the potential $V^{l+\frac{1}{2}}$ and the wave function $\Psi^l_p$.
\item[$S2.$]Numerical integration: for $j=0,...,J$
\begin{equation}\label{eqID}
\displaystyle n^{l+1}_j= \sum_{p=0}^{P-1}\left(g(k_p)\left|\Psi^{l+1}_{p,j}\right|^2+g(k_{p+1})\left|\Psi^{l+1}_{p+1,j}\right|^2\right)\frac{\Delta k}{2}\,.
\end{equation}
\end{itemize} 
\item[]{\bf Coupling, Poisson equation}\\
Computation of the potential $V^{l+1}$ from the density $n^{l+1}$ by solving the Poisson equation \eqref{eqPot}:
\begin{equation}\label{potLp1}
\displaystyle \left\{
\begin{array}{ll} -\partial_x^2V^{l+1} = \frac{q^2}{\varepsilon}(n^{l+1}-n_D)\,, & x\in(0,L)\\[1.65mm]
V^{l+1}(0)=V^{l+1}(L)=0
\end{array}\right.\,.
\end{equation}
\end{itemize} 
\underline{End Do}.
\end{algorithm}
In equation \eqref{eqID}, $n^{l+1}_j$ denotes the approximation of the density $n^{l+1}(x_j)$ and $\Psi^{l+1}_{p,j}$ the approximation of the wave function $\Psi^{l+1}_p(x_j)$. The linear equation \eqref{potLp1} is discretized with finite difference and solved by a simple matrix inversion. As in the stationary case, the steps $S1$ and $S2$ to compute the density will be crucial for the treatment of the resonance peaks.\\
In the algorithm proposed in \cite{Pi}, that we will call Direct Resolution, the step $S1$ is performed as follows: the function $\Psi^{l+1}_p$ is computed on $[0,L]$ using $\Psi^{l}_p$ and $V^{l+\frac{1}{2}}$ by the resolution of \eqref{schInst} with the Crank-Nicolson scheme \eqref{CranNich}, where the potential $Q$ is equal to $U+V$. Since the support of the initial data $\Phi_{k_p}$ is not included in $(0,L)$, the discrete transparent boundary conditions for \eqref{CranNich} are the non-homogeneous ones \eqref{dtbcNHg}\eqref{dtbcNHd} or their simplified version \eqref{dtbcNHg_cut}\eqref{dtbcNHd_cut}. In the Direct Resolution, no particular treatment of the resonances is performed in the step $S1$ and the accuracy of the method is provided in the step $S2$ by the imposition of a small uniform grid spacing $\Delta k$ in the trapezoidal rule \eqref{eqID}. Indeed, if the grid spacing is small only near the initial resonant frequencies, the refined mesh will loose its relevance because of the time evolution of the resonance. Therefore, the stationary Adaptive Method does not extend to the transient regime and, if no additional treatment is performed (a strategy would be to choose a small frequency grid spacing only in the regions where will live the resonance peaks), the number of Schr\"odinger equations to be solved have to be important in the time-dependent Direct Resolution. Moreover, in addition to a matrix inversion, which is equally required in the stationary case, each time-dependent Schr\"odinger equation requires discrete convolutions in the boundary conditions. For the exact boundary conditions \eqref{dtbcNHg}\eqref{dtbcNHd}, the number of operations required by these convolutions is proportional to the time iteration number $l$. For the simplified boundary conditions \eqref{dtbcNHg_cut}\eqref{dtbcNHd_cut}, the size of the convolutions can be reduced but it can not be too small in order to preserve the accuracy of the method. This increases the numerical cost of each time step and makes the time-dependent simulations very long: actually, due to the condition \eqref{eqstabt_toutk}, the time step must be chosen much smaller than the scaling time of the system and the number $l$ can become very big. In this context, it is important to look for an adapted treatment of the resonance peaks to reduce the number of frequency points.\\
In the following section, we propose a One Mode Approximation method which extends the method proposed in section \ref{sec_1mSt} to the time-dependent case.

\subsection{The time-dependent One Mode Approximation}\label{sec_1mInst}

Let us start with the description of the decomposition of a wave function $\Psi_k(t)$ solution to the transient Schr\"odinger equation \eqref{schInst} for a given frequency $k$. We look for $\Psi_k(t)$ of the form:
\begin{equation}\label{form1mIns1}
\displaystyle \Psi_k(t)=\Psi^{nr}_k(t)+\Psi^r_k(t)\,,
\end{equation}
where the non resonant part $\Psi^{nr}_k(t)$ solves the transient Schr\"odinger equation:
\begin{equation}\label{SchroInstNR}
\displaystyle \left\{
\begin{array}{ll}i\hbar\partial_t\Psi^{nr}_k(t) = [-\frac{\hbar^2}{2m}\partial_x^2+U_{fill}(t)+V(t)]\Psi^{nr}_k(t)\,, & x\in\R\\[1.65mm]
\Psi^{nr}_k(0) = \Phi^{nr}_k
\end{array}\right.\,,
\end{equation}
with filled potential $U_{fill}(t)=U(t)+v_0{\bf1}_{[a_3,b_3]}$ and $\Phi^{nr}_k$  is solution to \eqref{SchroStatNR}. Using the same argument as in the stationary regime, we suppose that $\Psi^r_k(t)$ is proportional to the resonant mode $u(t)$ corresponding to the first resonance $z(t)$ of the Hamiltonian $H(t)$ at time $t$: in other words
\begin{equation}\label{form1mIns2}
\displaystyle \Psi^r_k(t,x) =\lambda_k(t)u(t,x)\,, \quad x\in[0,L]\,,
\end{equation}
where 
\begin{equation}\label{eqRes}
\displaystyle [-\frac{\hbar^2}{2m}\partial_x^2+U(t)+V(t)]u(t)=z(t)u(t) \quad \textrm{ and } \quad \int_0^L|u(t,x)|^2dx=1\,. 
\end{equation}
If $u_I$ and $\theta_k$ are respectively the resonant mode and the proportionality coefficient given by \eqref{eqResI} and \eqref{eqTheta}, then it holds from \eqref{form1mStat1}\eqref{form1mStat2}:
\begin{equation}\label{form1mStatc}
\displaystyle \Phi_k(x)=\Phi^{nr}_k(x)+\theta_ku_I(x)\,, \quad x\in[0,L]\,.
\end{equation}
Comparing \eqref{form1mIns1}\eqref{form1mIns2} and \eqref{form1mStatc}, the initial condition in \eqref{schInst} and \eqref{SchroInstNR} imply $\lambda_k(0)=\theta_k$. Moreover, injecting \eqref{form1mIns1}\eqref{form1mIns2} in the transient Schr\"odinger equation \eqref{schInst}, and utilizing the equations \eqref{SchroInstNR} and \eqref{eqRes}, we get the following equation on $\lambda_k(t)$:
\[
\displaystyle [i\hbar\lambda_k'(t)-z(t)\lambda_k(t)]u(t,x)+i\hbar\lambda_k(t)\partial_tu(t,x)=-v_0{\bf 1}_{[a_3,b_3]}(x)\Psi^{nr}_k(t,x)\,.
\]
Multiplying the previous equation by $\overline{u}(t,x)$ and integrating on $(0,L)$, it follows:
\begin{equation}\label{edoLa}
\displaystyle \left\{
\begin{array}{l}
\lambda_k'(t) + [\frac{i}{\hbar}z(t) +
  \int_0^L\partial_tu(t,x)\overline{u}(t,x)dx]\lambda_k(t) = \frac{i}{\hbar}v_0\int_{a_3}^{b_3}\Psi_k^{nr}(t,x)\overline{u}(t,x)dx\\[1.65mm]
\lambda_k(0)=\theta_k
\end{array}\right.\,.
\end{equation}
\begin{remark}\label{rq_picla}
In the linear case, obtained by taking $V(t)=0$ in \eqref{schInst}, it holds for $t>0$: $u(t)=u_\infty$ and $z(t)=E_\infty-i\Gamma_\infty$, where $E_\infty\in\R$ and $\Gamma_\infty>0$. Then, under the assumption $\Psi_k^{nr}(t,x)=\tilde{\Psi}_k^{nr}(x)e^{-i\frac{\varepsilon_k^\infty}{\hbar}t}$, where
\begin{equation}\label{eqDispt}
\displaystyle \varepsilon^\infty_k=\left\{\begin{array}{ll} \frac{\hbar^2k^2}{2m}\,, & k\geq 0\\[1.5mm]
\frac{\hbar^2k^2}{2m}-B_{\infty}\,, & k<0\end{array}\right.\,,
\end{equation}
equation \eqref{edoLa} is an ODE which homogeneous solution oscillates at the angular frequency $\frac{E_\infty}{\hbar}$ and which source term oscillates at the angular frequency $\frac{\varepsilon_k^\infty}{\hbar}$. Therefore, $\lambda_k(t)$ has a peak at the frequencies $k$ such that $\varepsilon_k^\infty=E_\infty$. Due to the initial condition in \eqref{edoLa}, one expects a coexistence of this peak with a second peak at the frequencies $k$ such that $E_k=E_I$ with decay rate $\frac{\Gamma_\infty}{\hbar}$, where $E_k$ and $E_I$ are given by \eqref{eqDisp} and \eqref{eqZI}. This is verified numerically in section \ref{sec2pics}.
\end{remark}

\subsubsection{Step $S1$}

The aim here is to compute the elements of the decomposition \eqref{form1mIns1}\eqref{form1mIns2} at the time $t=t^{l+1}$. At this point of the algorithm, we do not have the value of the potential $V^{l+1}$ and the corresponding resonant mode $u^{l+1}$ can not be computed. Therefore, we will make the following approximation for the wave function at time $t^{l+1}$ and frequency $k$:
\begin{equation}\label{dec_psiIns}
\displaystyle \Psi_k^{l+1}=\Psi^{nr,l+1}_k+\lambda_k^{l+1}u^{l+\frac{1}{2}}\,.
\end{equation}
Then the step $S1$ writes as follows: suppose that the quantities $\Psi^{nr,l}_k$ and $\lambda_k^{l}$ are known (at $l=0$, it is given by the initial decomposition \eqref{form1mStatc}).\\
As it is done in the Direct Resolution for $\Psi_k^{l+1}$, the non resonant wave function $\Psi^{nr,l+1}_k$ at time $t^{l+1}$ and frequency $k$ is computed on $[0,L]$ using $\Psi^{nr,l}_k$ and $V^{l+\frac{1}{2}}$ by the resolution of \eqref{SchroInstNR} with the Crank-Nicolson scheme \eqref{CranNich} where the potential $Q$ is equal to $U_{fill}+V$ and the initial data $\Phi$ is equal to $\Phi^{nr}_k$. Equation \eqref{CranNich} comes with the transparent boundary conditions \eqref{dtbcNHg}\eqref{dtbcNHd} or \eqref{dtbcNHg_cut}\eqref{dtbcNHd_cut}.\\
The resonant mode $u^{l+\frac{1}{2}}$ and the resonance $z^{l+\frac{1}{2}}$ are computed using the method presented in section \ref{sec_calres} where the potential $Q$ is equal to $U^{l+\frac{1}{2}}+V^{l+\frac{1}{2}}$ and $U^{l+\frac{1}{2}}=U(t^{l+\frac{1}{2}})$. The initial guess for $(u^{l+{1\over 2}}, z^{l+{1\over 2}})$ are the first eigenfunction and eigenenergy  of the Hamiltonian 
\begin{equation}\label{eq_Hamd}
\displaystyle [-\frac{\hbar^2}{2m}\partial_x^2+U^{l+\frac{1}{2}}+V^{l+\frac{1}{2}}]
\end{equation}
with homogeneous Dirichlet boundary conditions at $a_2$ and $b_2$.\\
To achieve the step $S1$, we have to compute the coefficient $\lambda_k^{l+1}$ at time $t^{l+1}$ and frequency $k$. It is obtained from $\lambda_k^{l}$, $\Psi^{nr,l}_k$, $\Psi^{nr,l+1}_k$, $z^{l+\frac{1}{2}}$ and $u^{l+\frac{1}{2}}$ by the resolution of equation \eqref{edoLa} with a trapezoidal rule. This leads to the following iteration:
\begin{multline}\label{trapTh}
\displaystyle (\lambda^{l+1}_k-\lambda^l_k)/\Delta t + [\frac{i}{\hbar}z^{l+\frac{1}{2}} + \int_0^L\partial_tu^{l+\frac{1}{2}}\overline{u^{l+\frac{1}{2}}}dx](\lambda^{l}_k+\lambda^{l+1}_k)/2\\ = \frac{iv_0}{2\hbar}\int_{a_3}^{b_3}\left[\Psi^{nr,l}_k+\Psi^{nr,l+1}_k\right]\overline{u^{l+\frac{1}{2}}}dx\,,
\end{multline}
where $\partial_tu^{l+\frac{1}{2}}$ stands for $\partial_tu(t^{l+\frac{1}{2}})$. By adequately fixing the resonant mode phase, the quantity
\begin{equation}\label{eqmu}
\displaystyle \mu^{l+1/2} := \int_0^L \partial_tu^{l+1/2}\overline{u^{l+1/2}} dx
\end{equation}
appearing in \eqref{trapTh} is fitted to zero. Given $\tilde{u}^{l+1/2}$ resonant mode of lower energy of the Hamiltonian \eqref{eq_Hamd} and verifying $\Vert\tilde{u}^{l+1/2}\Vert_{L^2(0,L)}=1$, we look for $u^{l+1/2}$ of the form $u^{l+1/2}=\tilde{u}^{l+1/2}e^{i\varphi^{l+1/2}}$ where $\varphi^{l+1/2} \in \R$ is such that $\mu^{l+1/2}$ can be neglected. We note first that we have the approximation
\[
\displaystyle \mu^{l+1/2} = \frac{1}{2\Delta t} \int_0^L(u^{l+1/2}-u^{l-1/2})\overline{(u^{l+1/2}+u^{l-1/2})}dx\,,
\]
which becomes
\begin{equation}\label{approw_mu}
\displaystyle \mu^{l+1/2} =  \frac{i}{\Delta t}\textrm{Im}[\int_0^Lu^{l+1/2}\overline{u^{l-1/2}}dx]
\end{equation}
under the condition $\Vert u^{l-1/2} \Vert_{L^2(0,L)}=\Vert u^{l+1/2} \Vert_{L^2(0,L)}=1$. Then, defining 
\begin{equation}\label{eq_omd}
\displaystyle \omega^{l+1/2} = \int_0^L\tilde{u}^{l+1/2}\overline{u^{l-1/2}}dx\,,
\end{equation}
 we choose
\[
\displaystyle e^{i\varphi^{l+1/2}} = \frac{\overline{\omega^{l+1/2}}}{|\omega^{l+1/2}|}
\]
and it follows:
\[ 
\displaystyle \int_0^Lu^{l+1/2}\overline{u^{l-1/2}}dx = e^{i\varphi^{l+1/2}}\int_0^L\tilde{u}^{l+1/2}\overline{u^{l-1/2}}dx = \omega^{l+1/2}e^{i\varphi^{l+1/2}} = |\omega^{l+1/2}| \in \R\,.
\]
As a consequence, $\textrm{Im}[\int_0^Lu^{l+1/2}\overline{u^{l-1/2}}dx] = 0$ and equation \eqref{approw_mu} shows that $u^{l+1/2}$ is such that $\mu^{l+1/2}$ is almost equal to $0$. In the applications considered in section \ref{sec_results}, the space integral in \eqref{eq_omd} is computed using a rectangle method.
\begin{remark}
As it is done in the stationary case, see Remark \ref{rqInitres}, the algorithm can be accelerated by choosing $(u^{l-{1\over 2}}, z^{l-{1\over 2}})$, for $l\geq 1$, as initial guess of the method to compute $(u^{l+{1\over 2}}, z^{l+{1\over 2}})$. When $l=0$, the initial guess is the resonant mode and the resonance given by the stationary solution. 
\end{remark}
\begin{remark}
Although it is harder to implement than an explicit Euler method, the semi-implicit trapezoidal rule \eqref{trapTh}, used to solve the stiff-like equation \eqref{edoLa}, is required for stability purposes. Indeed, the imaginary part of the resonance $z(t)$ being strictly negative, the ODE \eqref{edoLa} is similar to the test equation: $y'=\alpha y$ where $\textrm{Re}\,(\alpha) < 0$, for which the trapezoidal rule solution is bounded independently of the time step whereas the Euler explicit solution is conditionally bounded. In particular, as explained above, the quantity $\mu^{l+1/2}$ defined in \eqref{eqmu} can been fitted to zero in \eqref{trapTh}, and a straightforward calculation shows that the resulting algorithm is unconditionally stable, in the sens that for all $N\geq 0$
\[
\left\vert\lambda^l_k\right\vert \leq \left\vert\lambda^0_k\right\vert + T\max_{0\leq n\leq N-1}\left\vert\frac{iv_0}{2\hbar}\int_{a_3}^{b_3}\left[\Psi^{nr,n}_k+\Psi^{nr,n+1}_k\right]\overline{u^{n+\frac{1}{2}}}dx\right\vert\,, \quad 0\leq l\leq N\,,  
\]
where $T=N\Delta t$.
\end{remark}

\subsubsection{Step $S2$}\label{sec_S2I}

The step $S2$ is performed using the decomposition \eqref{dec_psiIns} to find an approximation of the density $n^{l+1}_j$ in \eqref{eqID} which provides an adapted treatment of the resonant peaks. Like in section \ref{sec_1mSt}, we make the approximation:
\begin{equation}\label{eq.decpsi_disc1}
\displaystyle \left|\Psi_k^{l+1}(x)\right|^2=\left|\Psi^{nr,l+1}_k(x)\right|^2+\left|\lambda_k^{l+1}\right|^2\left|u^{l+\frac{1}{2}}(x)\right|^2\,.
\end{equation}
The non resonant wave function $\Psi^{nr}_k$ is regular with respect to the frequency $k$ and the integral of the first term at the r.h.s. of equation \eqref{eq.decpsi_disc1} can be computed using a trapezoidal rule with coarse frequency mesh. As it is noticed in Remark \ref{rq_picla}, the coefficient $\lambda_k$ has sharp peaks at the resonant frequencies and the frequency mesh is taken thin for the integral of the second term at the r.h.s of equation \eqref{eq.decpsi_disc1}. More precisely, we consider the frequency meshes $\{\hat{k}_{p'},\,p'=0,...,P'\}$ and $\{k_p,\,p=0,...,P\}$ introduced in section \ref{sec1mapp} and which have different scales adjusted using the ratio $\nu=\frac{P}{P'}\in\mathbb{N}$. Then, if we note $\Psi^{nr,l}_{p',j}$ the approximation of the non resonant wave function $\Psi^{nr,l}_{k_{\nu p'}}(x_j)$, $\lambda_p^l$ the approximation of the coefficient $\lambda^l_{k_p}$ and $u_j^{l+\frac{1}{2}}$ the approximation of the resonant mode $u^{l+\frac{1}{2}}(x_j)$, the formula for the density is:
\begin{multline}\label{eq_dbscint}
n_{j}^{l+1}=\sum_{p'=0}^{P'-1}\left(g(\hat{k}_{p'})\left|\Psi^{nr,l+1}_{p',j}\right|^2+g(\hat{k}_{p'+1})\left|\Psi^{nr,l+1}_{p'+1,j}\right|^2\right)\frac{\Delta \hat{k}}{2} \\
+\sum_{p=0}^{P-1}\left(g(k_{p})\left|\lambda_p^{l+1}\right|^2+g(k_{p+1})\left|\lambda_{p+1}^{l+1}\right|^2\right)\frac{\Delta k}{2}\left|u_j^{l+\frac{1}{2}}\right|^2\,.
\end{multline}
The number of Schr\"odinger equations to solve is reduced, $P'$ equations instead of $P$ equations for the Direct Resolution, reducing this way the numerical cost. However, this reduction implies that we have only access to the functions $\Psi^{nr,l}_{p',j}$ for $p'=0,...,P'$ and the computation of the coefficient $\lambda_{p}^{l+1}$, for $p\neq \nu p'$, requires an interpolation of the non resonant wave function to evaluate the source term in \eqref{trapTh}. If we set $\varepsilon^{\infty}_{p}=\varepsilon^\infty_{k_p}$ where $\varepsilon^\infty_{k}$ is given by \eqref{eqDispt}, it follows from Remark \ref{rq_picla} that the coefficient $\lambda_{p}^{l+1}$ has a peak when the energy $\varepsilon^{\infty}_{p}$ is equal to $E^{l+1}$ (the lower resonant energy corresponding to the potential $V^{l+1}$) only if the approximation of the source term in \eqref{trapTh} includes the time oscillations of the non resonant wave function. In particular, a polynomial interpolation of the non resonant wave function is not adapted. Solving equation \eqref{trapTh} and performing an oscillatory interpolation of the non resonant wave function, we obtain the following suitable algorithm which extends the formula \eqref{eqThetaIp} to the time-dependent case: for $p=0,...,P$
\begin{multline}\label{eqLaint}
\lambda^{l+1}_p=\frac{1}{1+i\frac{\Delta t z^{l+\frac{1}{2}}}{2\hbar}}\left[ \left(1-i\frac{\Delta t z^{l+\frac{1}{2}}}{2\hbar}\right)\lambda^{l}_p \right. \\ \left. + \frac{iv_0\Delta t}{2\hbar}\sum_{j\in w}\left(\tilde{\Psi}^{nr,l}_{p',j}e^{-i\omega^{\infty}_pt^l}+\tilde{\Psi}^{nr,l+1}_{p',j}e^{-i\omega^{\infty}_pt^{l+1}}\right)\overline{u^{l+\frac{1}{2}}_j}\Delta x\right]\,,
\end{multline}
where $w$ is defined in section \ref{sec1mapp} and $0\leq p'\leq P'$ is the integer such that $p=\nu p'+r$ for some $0\leq r\leq \nu-1$. The amplitude-like function $\displaystyle\tilde{\Psi}^{nr,l}_{p',j}$ appearing in equation \eqref{eqLaint} is defined by
\begin{equation*}
\tilde{\Psi}^{nr,l}_{p',j} = \Psi^{nr,l}_{p',j}e^{i\omega^{\infty}_{\nu p'}t^l}
\end{equation*}
and
\begin{equation}\label{eqwr}
\omega^{\infty}_p = \frac{\varepsilon^{\infty}_p}{\hbar}\,.
\end{equation}
\begin{remark}
The condition \eqref{eqstabt_toutk} insures that the time step is small enough with respect to one period of the oscillating term $e^{-i\omega^{\infty}_pt}$ in \eqref{eqLaint}. The time-dependent One Mode Approximation extends the method in section \ref{sec1mapp} and it is numerically verified that it is stable in time and allows long time simulations, see section \ref{sec_compt}. On the other hand, due to the function $e^{-i\omega^{\infty}_pt}$, a time-dependent extension of the method presented in section \ref{sec1mInt} is not obvious and gives a non-stable algorithm if no particular treatment is performed.
\end{remark}
\begin{remark}
In section \ref{sec-discangf}, it is noticed that the accuracy of the algorithm to solve the time-dependent Schr\"odinger equation is improved by replacing the angular frequencies \eqref{eqwg}\eqref{eqwd} by the discrete ones \eqref{eqwgd}\eqref{eqwdd} in the boundary conditions \eqref{dtbcNHg}\eqref{dtbcNHd}. Similarly, the accuracy of the time-dependent One Mode Approximation is improved when the angular frequency \eqref{eqwr} is replaced by the discrete one below:
\begin{equation}\label{eqwrd}
\omega^{\infty}_p = \frac{2}{\Delta t}\arctan\left(\frac{\Delta t \varepsilon^{\infty}_p}{2\hbar}\right)\,.
\end{equation}
\end{remark}

\section{Results}\label{sec_results}

The physical parameters used for the numerical computations are gathered in the following array:
\begin{table}[!ht]\begin{center}
\begin{tabular*}{0.48\textwidth}{@{\extracolsep{\fill}} ll | ll}   
   Rel. el. mass & $0.067$       &  Rel. permitivity & $11.44$ \\
   Temperature & $300\, K$ & Fermi level $E_F$ & $6,7097\times10^{-21}\,J$\\ 
   Donor density $n_D^1$ & $10^{24}\, m^{-3}$\\
Donor density $n_D^2$ & $5\times10^{21}\, m^{-3}$
\end{tabular*}
\end{center}\end{table}\\
For all the tests, the two barriers have the same size which is equal to the size of the well. The external potential is determined by Figure \ref{figVe} and the data below:
\begin{center}
\begin{tabular}{c c c c c c c c}
\hline 
$L\,(nm)$   & $a_1\,(nm)$  & $a_2\,(nm)$  & $a_3\,(nm)$	 & $b_3\,(nm)$  & $b_2\,(nm)$  & $b_1\,(nm)$  & $v_0\,(eV)$\\
\hline 
$135$ &	$50$ & $60$	& $65$ & $70$ & $75$ & $85$ & $0.3$\\
\hline
\end{tabular}
\end{center}
In all the simulations, we fixed $\kappa=\sqrt{\frac{2m}{\hbar^2}(E_F+7k_BT)}$ and the number of space points is equal to $J=300$ which is such that the condition \eqref{eqstabx_toutk} is verified.

\subsection{Computation of resonances}\label{sec_compres}

\begin{figure}
\begin{center}
\includegraphics[width=0.8\linewidth]{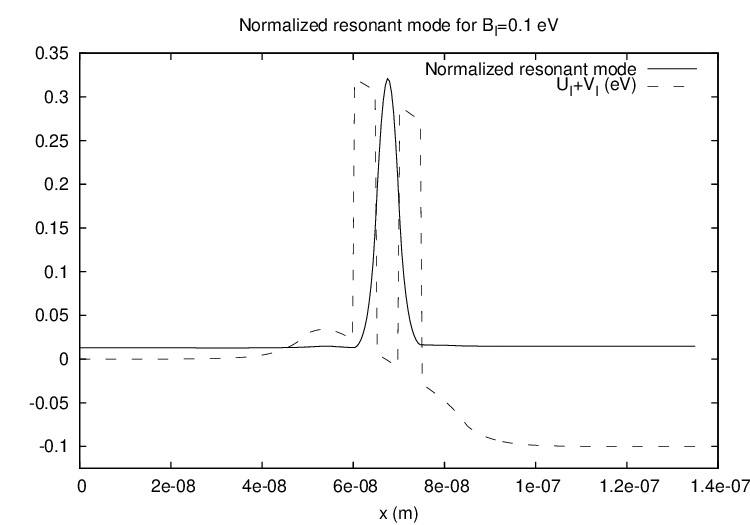}
\caption{Representation of the potential $U_I+V_I$ (dashed line) and the corresponding normalized resonant mode $\frac{|u_I(x)|}{|u_I|_2}$ (full line) for $B_I=0.1\;eV$.}\label{figmodReso}
\end{center}
\end{figure}

In this section, we give the numerical values of the resonance of lower energy obtained using the algorithm presented in section \ref{sec_calres} for different biases $B_I$. The potential $Q$ is equal to $U_I+V_I$ where $U_I$ is the external potential corresponding to $B_I$ and $V_I$ the corresponding solution to the Schr\"odinger-Poisson system \eqref{VpGe}\eqref{eqPotst}\eqref{denStat} given by the Direct Resolution, with $P=3750$, as described in section \ref{secComp}. We obtain the following results:
\begin{center}
\begin{tabular}{|c|c|c|c|c|}
\hline 
$B_I\,(eV)$   & $N_{cv}$  & $E_D\,(meV)$ & $E_I\,(meV)$ & $\Gamma_I/E_I$\\
\hline 
$0$       & $3$       & $135.14$   & $134.71$   & $1.7030\times 10^{-3}$ \\
\hline
$0.1$ & $3$	& $88.524$ & $88.090$ & $2.7547\times 10^{-3}$ \\
\hline
\end{tabular}
\end{center}
where the number of iterations before convergence $N_{cv}$ denotes the iteration $n$ of the algorithm such that $\left|M(z^n)u^n\right|_2<10^{-15}$, $E_D$ denotes the fundamental energy of the Dirichlet Hamiltonian \eqref{hamBar} and the resonance is equal to $z_I=E_I-i\Gamma_I$. The modulus of the normalized resonant mode for $B_I=0.1\,eV$ is represented in Figure \ref{figmodReso}.\\ 
The present section provides a numerical verification of the following theoretical results, which were used to write the One Mode Approximation: $E_D$ is close to $E_I$, $\Gamma_I$ is small with respect to $E_I$ and the restriction of the resonant mode to the interval $[0,L]$ is small outside the island $[a_2,b_2]$.

\subsection{The stationary regime}\label{secComp}

\begin{figure}
\begin{center}
\includegraphics[width=0.8\linewidth]{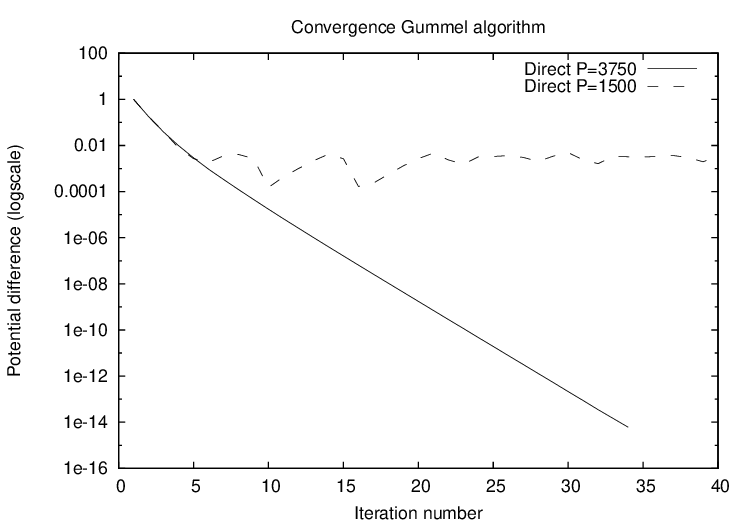}
\caption{Potential difference $\Delta$ defined by \eqref{eq_erreur} with respect to the number of iterations of Algorithm \ref{algGum} for the Direct Resolution for $B_I=0\,eV$ and for the values $P=1500$ and $P=3750$.}\label{cv_dirdec}
\end{center}
\end{figure}

\begin{figure}
\begin{center}
\begin{tabular}{ll}
$B_I = 0\,eV$ & \\
\includegraphics[width=0.48\linewidth]{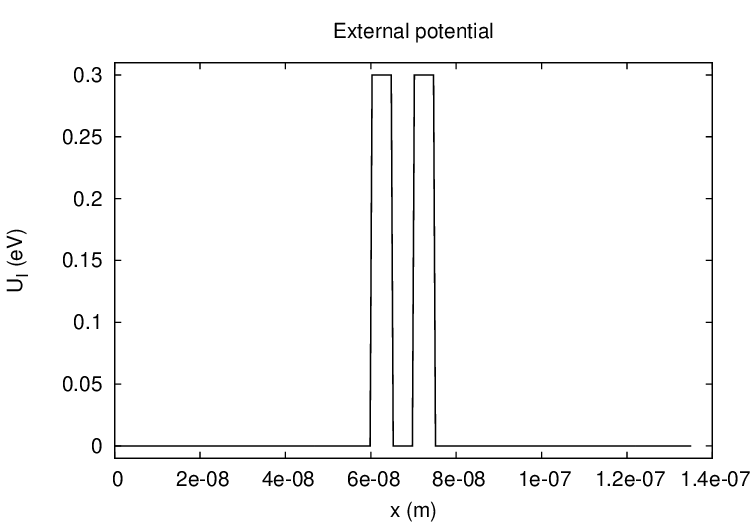}&
\includegraphics[width=0.48\linewidth]{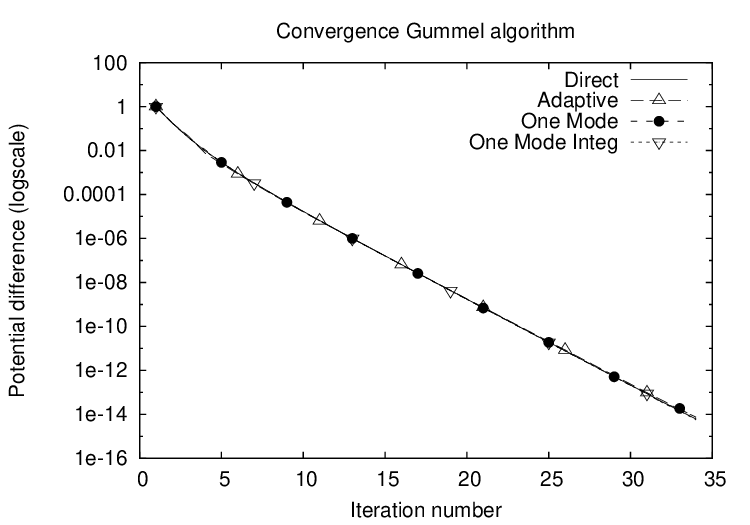}\\
\includegraphics[width=0.48\linewidth]{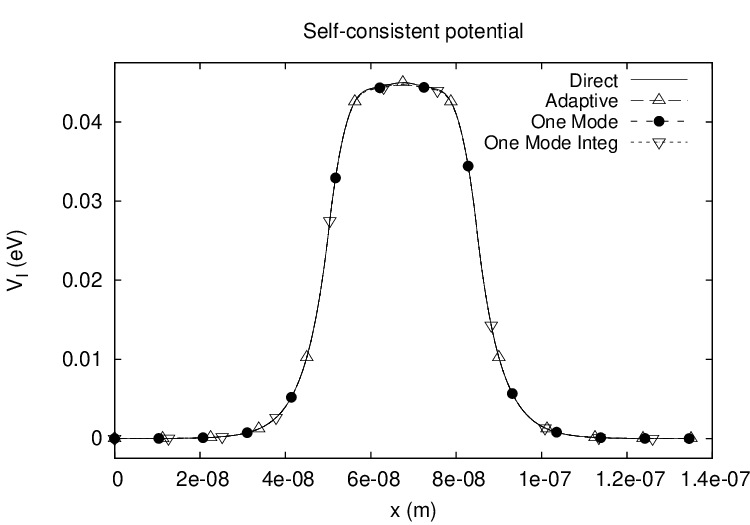}&\includegraphics[width=0.48\linewidth]{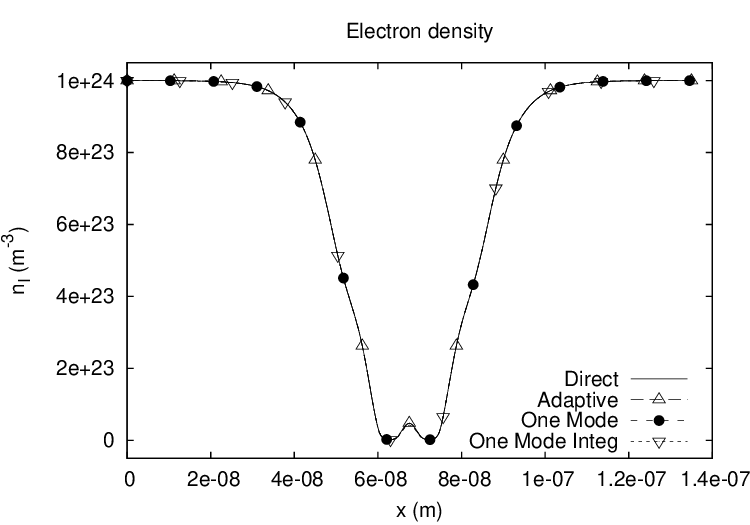}\\
$B_I = 0.1\,eV$ & \\
\includegraphics[width=0.48\linewidth]{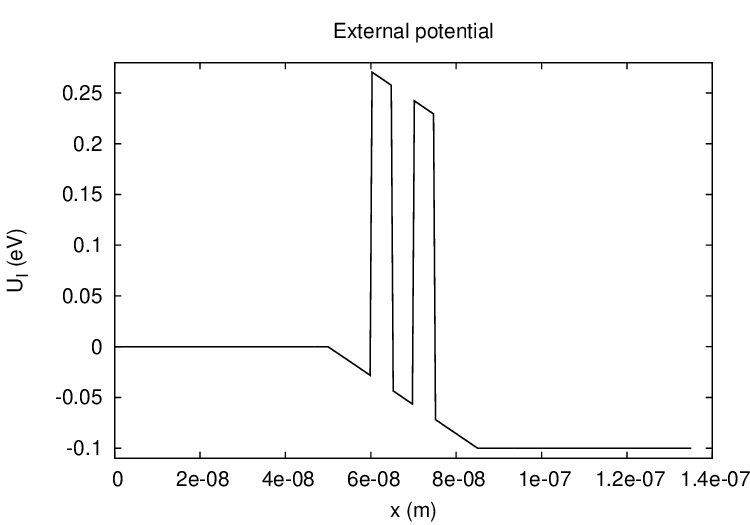}&
\includegraphics[width=0.48\linewidth]{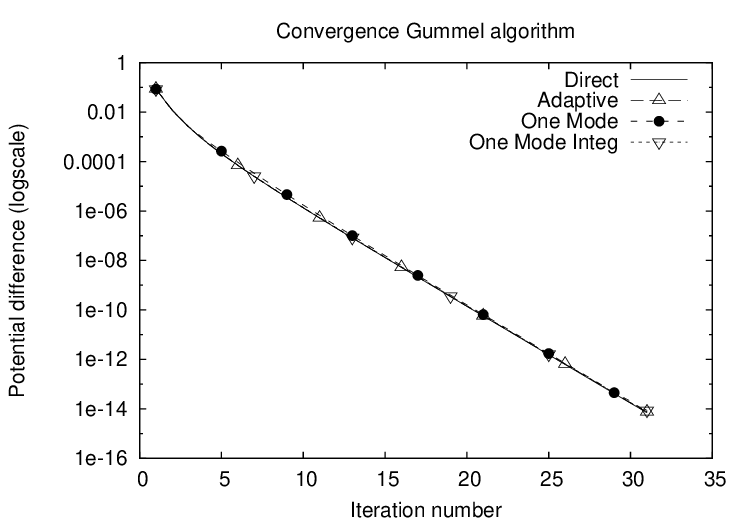}\\
\includegraphics[width=0.48\linewidth]{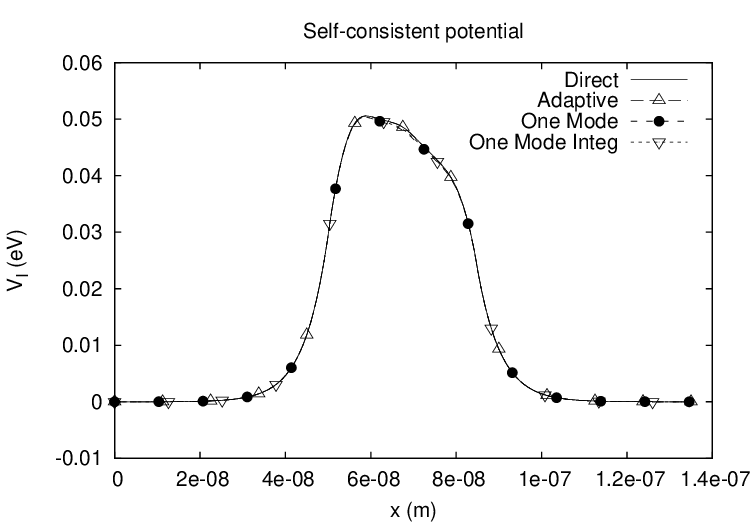}&\includegraphics[width=0.48\linewidth]{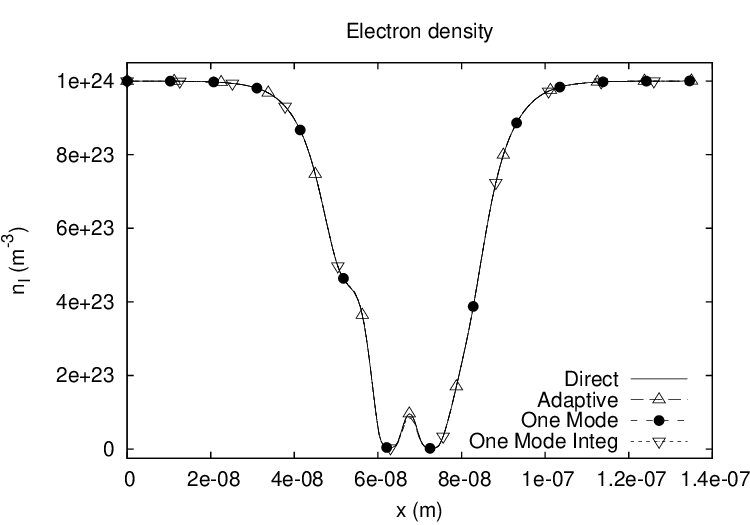}
\end{tabular}
\caption{For the biases $B_I = 0\,eV$ and $B_I = 0.1\,eV$ and, for the methods considered, are represented from left to right: the external potential $U_I$, the corresponding potential difference $\Delta$ defined by \eqref{eq_erreur} with respect to the number of iterations of Algorithm \ref{algGum}, the self-consistent potential $V_I$ and the electron density $n_I$.}\label{graph4meth}
\end{center}
\end{figure}

We give here a numerical comparison of the methods presented in section \ref{sec_stat}: the Direct Resolution, Adaptive Method, One Mode Approximation and One Mode Approximation Integ methods. The Direct Resolution plays the role of the reference method. For all the tests, the reference potential in \eqref{itGum} is given in electron-volts by $V_{ref}=\frac{T k_B}{q}$.
\begin{remark}
For the bias $B_I=0$, the Gummel algorithm \ref{algGum} can be initialized at the potential $V_I^0=0$. For $B_I>0$, such an initialization does not converge. In that case, the convergence is obtained by taking the initial potential equal to the solution corresponding to $B_I=0$.
\end{remark}
As it was underlined in section \ref{sec_algStat}, the choice of the frequency mesh determines the convergence and the numerical cost of the method. In particular, for the Direct Resolution, it is verified numerically that, $P\geq3700$ is required to insure the convergence of the method for the biases $B_I=0\,eV$ and $B_I=0.1\,eV$. This fact is illustrated in Figure \ref{cv_dirdec} for $B_I=0\,eV$.\\
Indeed, for $P=1500$, after $40$ iterations the method fails to converge for a CPU time equal to $1.8241\,s$. For $P=3750$, the method converges: the difference $\Delta$ becomes less than $tol=10^{-14}$ after $34$ iterations for a CPU time equal to $3.8562\,s$. For the Direct Resolution, it appears that a big number of frequency points, and therefore a relatively high numerical cost, is required to have convergence. This remark enhances the importance of the other methods.\\
The results given by the four methods, for two different values of the bias $B_I$, are given in Figure \ref{graph4meth} and Table \ref{tab4meth}. For the Adaptive Method, the small frequency step size is equal to $\Delta k = \frac{2\kappa}{4500}$ and the large step size is equal to $12\Delta k$. The frequency step size is equal to the small one when the logarithmic derivative of the transmission coefficient, see \cite{Pi,Fa} for its explicit formula, is greater than $1.3L$. The number of frequency points changes from an iteration to the other, however it stays around a fixed number which is $587$ for $B_I=0\,eV$ and $535$ for $B_I=0.1\,eV$. For the other methods, the frequency mesh is given by the integers $P=3750$ and $P'=150$.\\
\begin{table}
\begin{center}
\begin{tabular}{|c|c|c|c|c|}
\hline 
 &   & $N_{cv}$ & CPU($s$) & Dist to reference ($\%$) \\
\hline
 $B_I=0\,eV$ & Direct & $34$ & $3.8642$ & /\\
\cline{2-5} 
 & Adaptive & $34$ & $0.6040$ & $1.3859\times 10^{-2}$ \\
\cline{2-5} 
  & One Mode & $34$ & $0.6440$ & $0.2561$ \\
\cline{2-5} 
  & One Mode Integ& $34$ & $0.6040$ & $0.3700$ \\
\hline
$B_I=0.1\,eV$ & Direct & $31$ & $3.3402$ & / \\
\cline{2-5} 
& Adaptive & $31$ & $0.4800$ & $4.9157\times 10^{-3}$ \\
\cline{2-5} 
 & One Mode & $31$ & $0.5560$ & $0.7689$  \\
\cline{2-5} 
 & One Mode Integ& $31$ & $0.4880$ & $0.5126$  \\
\hline
\end{tabular}
\caption{Comparison of the Direct Resolution, Adaptive Method, One Mode Approximation and One Mode Approximation Integ methods for the resolution of the stationary Schr\"odinger-Poisson system \eqref{VpGe}\eqref{eqPotst}\eqref{denStat} for $B_I=0\,eV$ and $B_I=0.1\,eV$.}\label{tab4meth}
\end{center}
\end{table}  
In Table \ref{tab4meth}, the integer $N_{cv}$ is the total number of iterations of Algorithm \ref{algGum} for the tolerance $\textrm{tol}=10^{-14}$. The column CPU denotes the CPU time of the method and the column {\it Dist to reference}, the $l^2$-norm relative distance to the reference method potential. The last is given by: 
\[
\displaystyle 100\frac{\left|V^{N_{cv}}-V_{di}\right|_2}{\left|V_{di}\right|_2}\,,
\]
where $V_{di}$ is the potential of the Direct Resolution.\\
We notice that the alternative methods converge as fast as the reference method in term of number of iterations, and give a solution very close to the reference solution, for a smaller CPU time. For the alternative methods, by adequately fixing the parameters, and possibly allowing a bigger value in the column {\it Dist to reference} in Table \ref{tab4meth}, the frequency mesh can be chosen such that the CPU time is reduced. Anyway, for the frequency meshes chosen here, the alternative methods are an improvement to the reference method since their CPU times are lower than the CPU time of the Direct Resolution for the lowest value of $P$ such that the convergence of the method is ensured for the biases considered.\\
Although the One Mode Approximation and One Mode Approximation Integ methods require less frequency points than the Adaptive Method, the computation of the second term on the r.h.s. of equation \eqref{eq_densresst} increases the numerical cost of the two first methods, such that the three alternative methods have an equivalent CPU time. On the other hand, compared to the other alternative methods, the Adaptive Method provides a self-consistent potential which is closer to the Direct Resolution potential (this is due to the fact that the step $S1$ is performed in the same way for the two last methods). However, in the present work, the One Mode Approximation is more interesting than the other alternative methods since, as remarked in section \ref{sec_instat}, the last can not be generalized to the time-dependent case.

\subsection{The transient regime}

For all the tests, the external potential $U(t)$ is the time-dependent external potential corresponding to the bias $B(t)$ defined in \eqref{B_marche} where $B_I=0\,eV$ and $B_{\infty}=0.1\,eV$.

\subsubsection{Direct Resolution: Toward optimal parameters}\label{sec_parInst}
 
\begin{figure}
  \begin{minipage}[c]{.5\linewidth}
    \includegraphics[width=\linewidth]{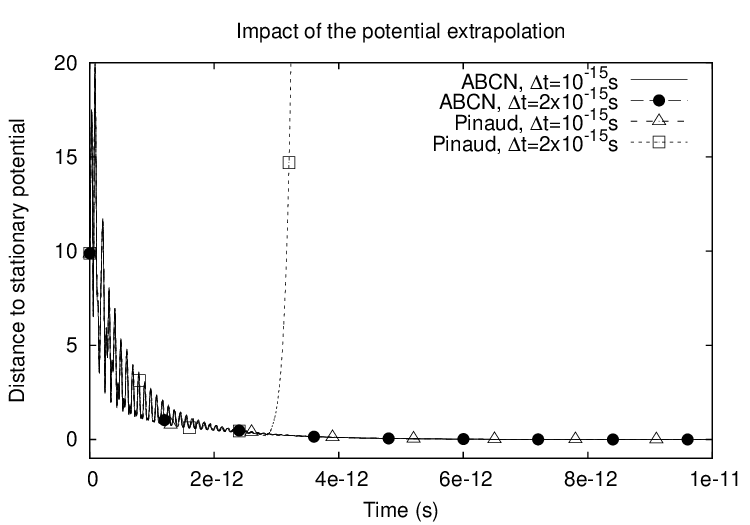}
  \end{minipage}
  \begin{minipage}[c]{.42\linewidth}
    {\footnotesize\begin{tabular}{|c|c|c|c|}
      \hline 
      & $l_{tot}$& CPU($s$) & Final dist ($\%$) \\
      \hline 
      ABCN, $\Delta t=10^{-15}\,s$ & $10^{4}$  & $3880.4$ & $4.4598\times10^{-4}$\\ 
      \hline
      ABCN, $\Delta t=2\times10^{-15}\,s$ & $5\times10^{3}$  & $1387.0$ & $6.0283\times10^{-4}$\\ 
      \hline
      Pinaud, $\Delta t=10^{-15}\,s$ & $10^{4}$ & $3864.8$ & $4.4609\times10^{-4}$ \\
      \hline
      Pinaud, $\Delta t=2\times10^{-15}\,s$ & $5\times10^{3}$ & $1440.7$ & $353.28$\\ 
      \hline
    \end{tabular}}
  \end{minipage}\\
  \begin{minipage}[c]{.5\linewidth}
    \includegraphics[width=\linewidth]{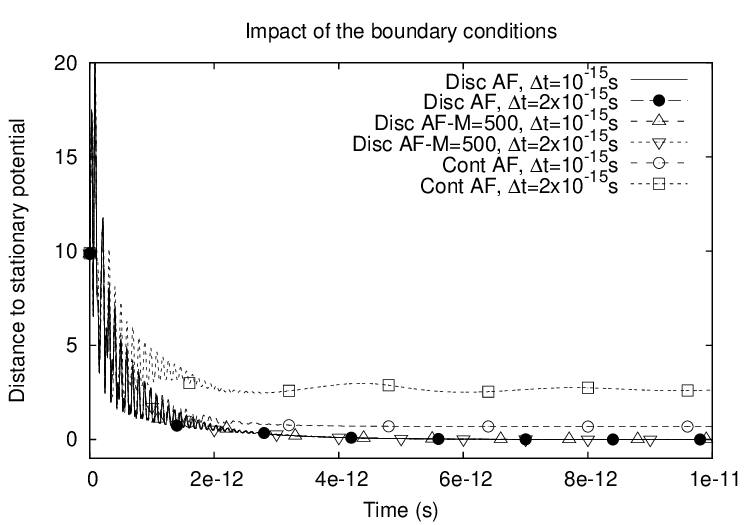}
  \end{minipage}
  \begin{minipage}[c]{.42\linewidth}
{\footnotesize\begin{tabular}{|c|c|c|c|}
      \hline 
      & $l_{tot}$& CPU($s$) & Final dist ($\%$) \\
      \hline 
      Disc AF, $\Delta t=10^{-15}\,s$ & $10^{4}$  & $3880.4$ & $4.4598\times10^{-4}$\\ 
      \hline
      Disc AF, $\Delta t=2\times10^{-15}\,s$ & $5\times10^{3}$  & $1387.0$ & $6.0283\times10^{-4}$\\
      \hline
      \specialcell{Disc AF-$M=500$,\\$\Delta t=10^{-15}\,s$} & $10^{4}$  & $1830.4$ & $8.0349\times10^{-3}$\\ 
      \hline
      \specialcell{Disc AF-$M=500$,\\$\Delta t=2\times10^{-15}\,s$} & $5\times10^{3}$  & $916.06$ & $4.5206\times10^{-3}$\\ 
      \hline
      Cont AF, $\Delta t=10^{-15}\,s$ & $10^{4}$ & $3867.8$ & $0.6874$ \\
      \hline
      Cont AF, $\Delta t=2\times10^{-15}\,s$ & $5\times10^{3}$ & $1388.1$ & $2.6214$\\ 
      \hline
    \end{tabular}}
  \end{minipage}\\
  \begin{minipage}[c]{.5\linewidth}
    \includegraphics[width=\linewidth]{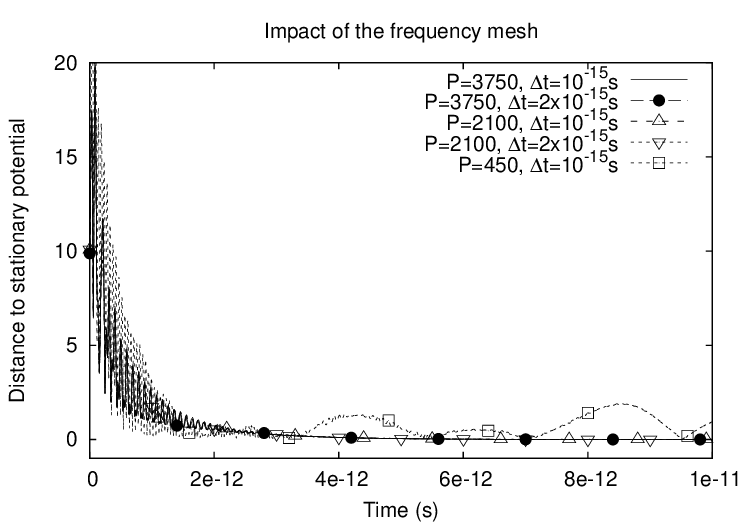}
  \end{minipage}
  \begin{minipage}[c]{.42\linewidth}
    {\footnotesize\begin{tabular}{|c|c|c|c|}
      \hline 
      & $l_{tot}$& CPU($s$) & Final dist ($\%$) \\
      \hline 
      $P=3750$, $\Delta t=10^{-15}\,s$ & $10^{4}$  & $3880.4$ & $4.4598\times10^{-4}$\\ 
      \hline
      $P=3750$, $\Delta t=2\times10^{-15}\,s$ & $5\times10^{3}$  & $1387.0$ & $6.0283\times10^{-4}$\\ 
      \hline
      $P=2100$, $\Delta t=10^{-15}\,s$ & $10^{4}$ & $2252.7$ & $9.7708\times10^{-3}$ \\
      \hline
      $P=2100$, $\Delta t=2\times10^{-15}\,s$ & $5\times10^{3}$ & $874.72$ & $4.5743\times10^{-3}$\\ 
      \hline
      $P=450$, $\Delta t=10^{-15}\,s$ & $10^{4}$ & $636.82$ & $0.94840$\\ 
      \hline
    \end{tabular}}
  \end{minipage}
  \caption{Evolution of the potential relative distance \eqref{eq_distv} to the stationary solution for the Direct Resolution. From top to bottom: Impact of the potential extrapolation, Impact of the boundary conditions and Impact of the frequency mesh.}\label{fig_dipar}
\end{figure}

The Direct Resolution presented in section \ref{sec_algInst} will play the role of the reference method and its CPU time is essential to measure the effective improvement provided by our alternative method. In this context, we show in the present section which choice of the parameters of the Direct Resolution provides the smallest CPU time while conserving a reasonable accuracy.\\
The initial data of the method is given by the stationary Direct Resolution for the initial bias $B_I=0\, eV$ and the parameters given in section \ref{secComp} where $P=3750$. The stationary potential corresponding to the final bias $B_{\infty}$ is computed in the same way and will be denoted $V_{\infty}$.\\
In Figure \ref{fig_dipar}, we depicted, for different parameters of the method, the time evolution of the $l^2$-norm relative distance to the potential $V_{\infty}$: 
\begin{equation}\label{eq_distv}
\displaystyle 100\frac{\left|V^{l}-V_{\infty}\right|_2}{\left|V_{\infty}\right|_2}\,,
\end{equation}
where $V^{l}$ is the Direct Resolution self-consistent potential corresponding to the iteration number $l$ of Algorithm \ref{AlgTra}. For each graphic of Figure \ref{fig_dipar}, is also given in the corresponding array: the total number of iterations $l_{tot}$, the CPU time of the method and the {\it Final dist} which corresponds to the value of the quantity \eqref{eq_distv} for $l=l_{tot}$.\\
In the first graphic of Figure \ref{fig_dipar}, the number of frequency points is $P=3750$ and we use the boundary conditions \eqref{dtbcNHg}\eqref{dtbcNHd} for the time-dependent Schr\"odinger equation where the angular frequencies are given by \eqref{eqwgd}\eqref{eqwdd}. When the extrapolation of the intermediary potential $V^{l+\frac{1}{2}}$  is given by \eqref{eqExtV_pi}, curve titled Pinaud, the Direct Resolution verifies the convergence criterion to the stationary potential $V_{\infty}$ only for the time step $\Delta t=10^{-15}s$. Although $\Delta t=2\times 10^{-15}s$ verifies time step limit \eqref{eqstabt_toutk} imposed by the linear Schr\"odinger equation, the method is clearly unstable for this time step and the loss of stability comes from \eqref{eqExtV_pi}. When the extrapolation is given by \eqref{eqExtV}, curve titled ABCN, the Direct Resolution potential convergences to $V_{\infty}$ for both time steps $\Delta t=10^{-15}s$ and $\Delta t=2\times 10^{-15}s$. Therefore, the extrapolation \eqref{eqExtV} provides a reduction of the CPU time. We remark however that the increase of the time step makes the amplitude of the transient oscillations of the potential slightly bigger.\\
In the second graphic of Figure \ref{fig_dipar}, the number of frequency points is $P=3750$ and we use the boundary conditions \eqref{dtbcNHg}\eqref{dtbcNHd} and the extrapolation \eqref{eqExtV} for the intermediary  potential $V^{l+\frac{1}{2}}$. When the angular frequencies are the discrete ones \eqref{eqwgd}\eqref{eqwdd}, curve titled Disc AF, the Direct Resolution potential convergences to $V_{\infty}$ and {\it Final dist} is small enough for both time steps $\Delta t=10^{-15}s$ and $\Delta t=2\times 10^{-15}s$. When the angular frequencies are the continuous ones \eqref{eqwg}\eqref{eqwd}, curve titled Cont AF, no instability is noticed, however, the value {\it Final dist} is much bigger than in the case Disc AF and increases significantly with the time step. In the case Disc AF, the computational cost is largely reduced while conserving a very satisfactory accuracy thanks to the simplified boundary conditions \eqref{dtbcNHg_cut}\eqref{dtbcNHd_cut} with $M=500$, curve titled Disc AF-$M=500$. Indeed, the CPU time is divided by more than two for $\Delta t=10^{-15}s$ and, for the two time steps, the value {\it Final dist} is less than $10^{-2}$.\\
In the third graphic of Figure \ref{fig_dipar}, we use the extrapolation \eqref{eqExtV} for the intermediary potential $V^{l+\frac{1}{2}}$ and the boundary conditions \eqref{dtbcNHg}\eqref{dtbcNHd} where the angular frequencies are given by \eqref{eqwgd}\eqref{eqwdd}. The different curves correspond to different values of the parameter $P$, curve titled $P=3750$, $P=2100$ and $P=450$. The case $P=2100$ illustrates the fact that a reasonable accuracy can be obtained even if the number of frequency points is under the limit $P=3700$ which insures the convergence of the stationary method. This fact is relatively stable with respect to the time step and allows a reduction of the CPU time of the method. However, the number of frequency points can not be chosen too small: in the case $P=450$, the value {\it Final dist} is bigger than $10^{-1}$ and the amplitude of the transient oscillations of the potential are large compared to the case $P=3750$. As in the stationary case, a big number of frequency points, and therefore a relatively high numerical cost, is required for the Direct Resolution, which enhances the importance of our alternative method.

\subsubsection{The One Mode Approximation}\label{sec_parInstr}

\begin{figure}
  \begin{minipage}[c]{.5\linewidth}
    \includegraphics[width=\linewidth]{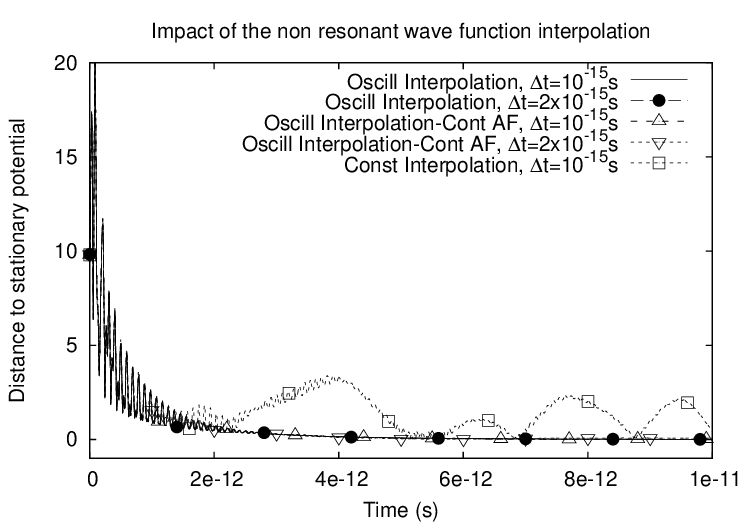}
  \end{minipage}
  \begin{minipage}[c]{.42\linewidth}
    {\footnotesize\begin{tabular}{|c|c|c|c|}
      \hline 
      & $l_{tot}$& CPU($s$) & Final dist ($\%$) \\
      \hline
      \specialcell{Oscill Interpolation,\\$\Delta t=10^{-15}\,s$} & $10^{4}$ & $496.03$ & $5.5306\times10^{-3}$\\ 
      \hline
      \specialcell{Oscill Interpolation,\\$\Delta t=2\times10^{-15}\,s$} & $5\times10^{3}$  & $204.29$ & $6.4993\times10^{-3}$\\ 
      \hline
      \specialcell{Oscill Interpolation-\\Cont AF, $\Delta t=10^{-15}\,s$} & $10^{4}$ & $497.31$ & $2.1769\times10^{-2}$ \\
      \hline
      \specialcell{Oscill Interpolation-\\Cont AF, $\Delta t=2\times10^{-15}\,s$} & $5\times10^{3}$  & $206.05$ & $7.7645\times10^{-2}$\\ 
      \hline
      \specialcell{Const Interpolation\\$\Delta t=10^{-15}\,s$} & $10^{4}$  & $506.45$ & $0.42549$\\ 
      \hline 
    \end{tabular}}
  \end{minipage}\\
  \begin{minipage}[c]{.5\linewidth}
    \includegraphics[width=\linewidth]{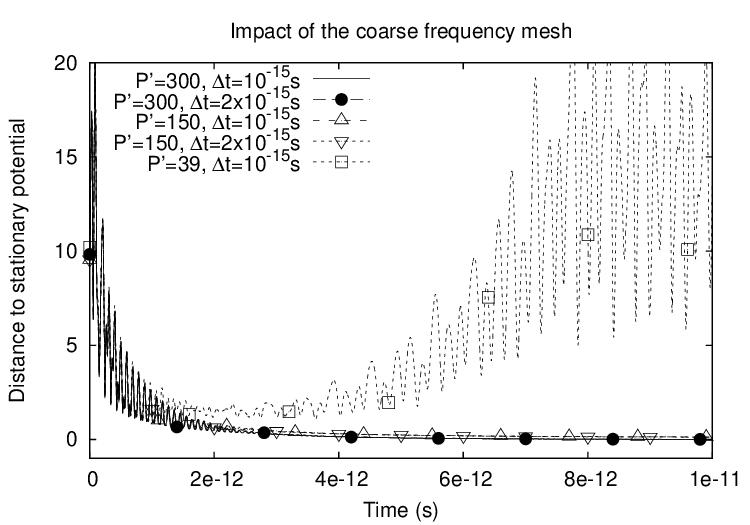}
  \end{minipage}
  \begin{minipage}[c]{.42\linewidth}
    {\footnotesize\begin{tabular}{|c|c|c|c|}
      \hline 
      & $l_{tot}$& CPU($s$) & Final dist ($\%$) \\
      \hline 
      $P'=300$, $\Delta t=10^{-15}\,s$ & $10^{4}$  & $496.03$ & $5.5306\times10^{-3}$\\ 
      \hline
      $P'=300$, $\Delta t=2\times10^{-15}\,s$ & $5\times10^{3}$  & $204.29$ & $6.4993\times10^{-3}$\\ 
      \hline
      $P'=150$, $\Delta t=10^{-15}\,s$ & $10^{4}$ & $344.06$ & $0.12238$ \\
      \hline
      $P'=150$, $\Delta t=2\times10^{-15}\,s$ & $5\times10^{3}$  & $152.50$ & $0.12444$\\ 
      \hline 
      $P'=39$, $\Delta t=10^{-15}\,s$ & $10^{4}$ & $242.87$ & $14.176$ \\
      \hline
    \end{tabular}}
  \end{minipage}
  \caption{Evolution of the potential relative distance \eqref{eq_distv} to the stationary solution for the One Mode Approximation. From top to bottom: Impact of the interpolation of the non resonant wave function and Impact of the coarse frequency mesh.}\label{fig_respar}
\end{figure}

Before comparing the two methods presented in section \ref{sec_instat}, we show how to adjust the parameters of the One Mode Approximation method, as we did in section \ref{sec_parInst} for the Direct Resolution.\\
The initial data of the method is given by the stationary One Mode Approximation with the initial bias $B_I=0\, eV$ and the parameters given in section \ref{secComp} where $P=3900$ and $P'=300$. The stationary potential $V_{\infty}$ corresponding to $B_{\infty}$ is computed in the same way.\\
For different parameters of the method, we depicted in Figure \ref{fig_respar} the time evolution of the $l^2$-norm relative distance \eqref{eq_distv} to the potential $V_{\infty}$ where $V^{l}$ is the One Mode Approximation self-consistent potential corresponding to the iteration number $l$ of Algorithm \ref{AlgTra}. For each graphic of Figure \ref{fig_respar}, we give in the corresponding array: the total number of iterations $l_{tot}$, the CPU time of the method and the number {\it Final dist} defined in section \ref{sec_parInst}. For all the curves represented, the time-dependent Schr\"odinger equation boundary conditions are given by \eqref{dtbcNHg}\eqref{dtbcNHd}, the angular frequencies by \eqref{eqwgd}\eqref{eqwdd} and the extrapolation of the intermediary potential by \eqref{eqExtV}. The number of frequency points corresponding to the thin mesh of our double scale algorithm is equal to $P=3900$.\\
In the first graphic of Figure \ref{fig_respar}, the number of frequency points corresponding to the coarse mesh is equal to $P'=300$. When the angular frequency $\omega^{\infty}_p$ in \eqref{eqLaint} is the discrete one \eqref{eqwrd}, curve titled Oscill Interpolation, the One Mode Approximation verifies the convergence criterion to the stationary potential $V_{\infty}$ for both time steps $\Delta t=10^{-15}s$ and $\Delta t=2\times10^{-15}s$. This is also true when the angular frequency $\omega^{\infty}_p$ is the continuous one \eqref{eqwr}, curve titled Oscill Interpolation-Cont AF. However, it is clear from the column {\it Final dist} that the discrete angular frequency provides a better convergence to the stationary potential than the continuous one. The constant interpolation of the non resonant wave function, curve titled Const Interpolation, is obtained by taking the angular frequency $\omega^{\infty}_p$ equal to $0$ in \eqref{eqLaint}. In that case, the One Mode Approximation is not accurate: the value {\it Final dist} is bigger than $10^{-1}$ and the potential oscillates when the time is bigger than $2\times10^{-12}s$.\\
In the second graphic of Figure \ref{fig_respar}, the angular frequency $\omega^{\infty}_p$ is given by \eqref{eqwrd}. The different curves correspond to different values of the parameter $P'$, curve titled $P'=300$, $P'=150$ and $P'=39$. The case $P'=150$ illustrates the fact that, if we allow a lower accuracy, the CPU time can be reduced by reducing the number $P'$ of Schr\"odinger equations to be solved. This fact is relatively stable with respect to the time step. However, it is clear from the case $P'=39$ that $P'$ can not be too small.\\
As in the Direct Resolution, the extrapolation equation \eqref{eqExtV} and the discrete angular frequencies \eqref{eqwgd}\eqref{eqwdd} present an improvement compared to the extrapolation equation \eqref{eqExtV_pi} and the angular frequencies \eqref{eqwg}\eqref{eqwd}. Here also, the CPU time can be reduced within a reasonable range of accuracy by using the simplified boundary conditions \eqref{dtbcNHg_cut}\eqref{dtbcNHd_cut} or by taking the number of frequency points, of the thin mesh, below the limit $P=3700$.

\subsubsection{Comparison for a suitable set of parameters}\label{sec_compt}

\begin{figure}
\begin{center}
\begin{tabular}{ll}
\includegraphics[width=0.49\linewidth]{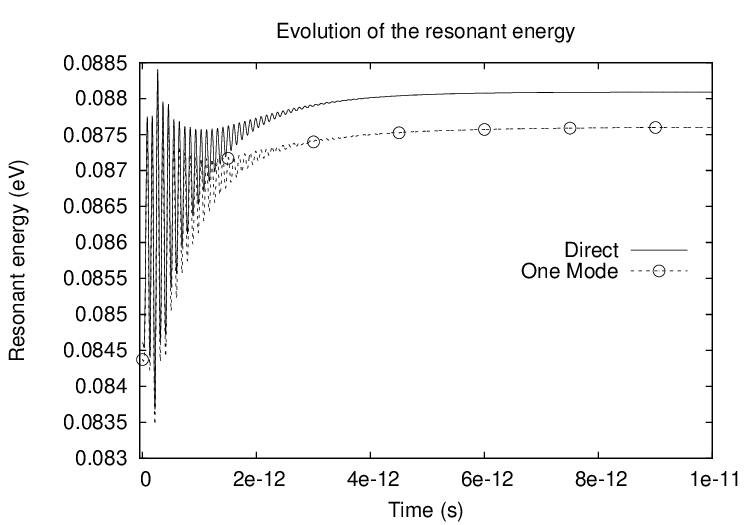}&
\includegraphics[width=0.49\linewidth]{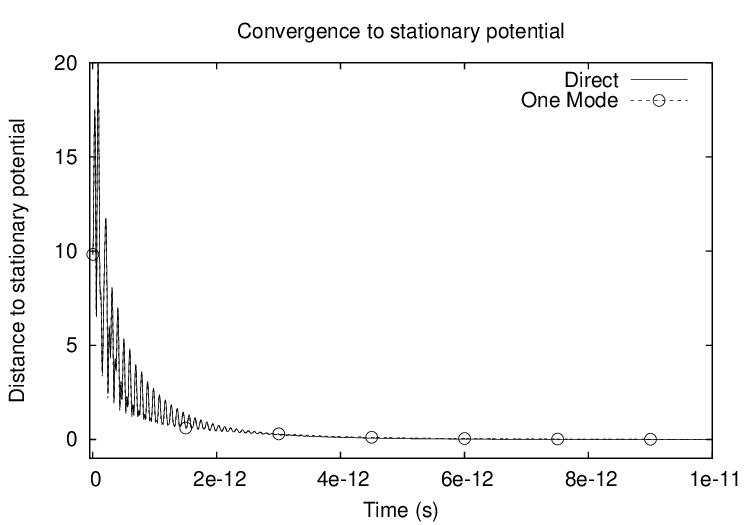}\\
\includegraphics[width=0.49\linewidth]{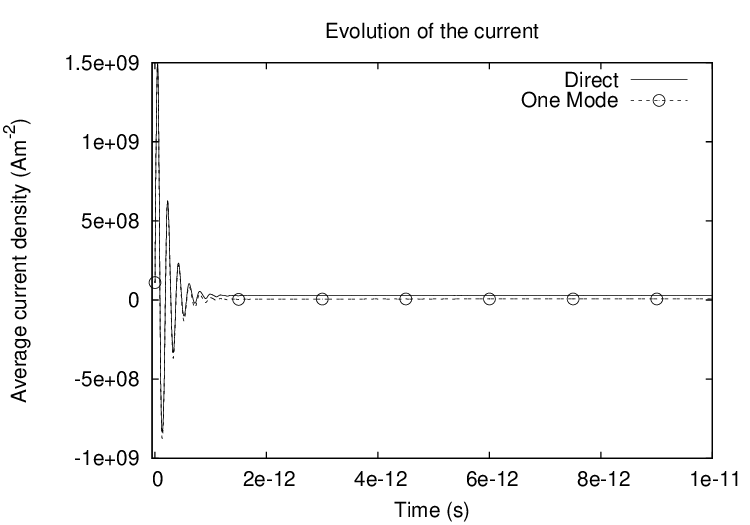}&
\includegraphics[width=0.49\linewidth]{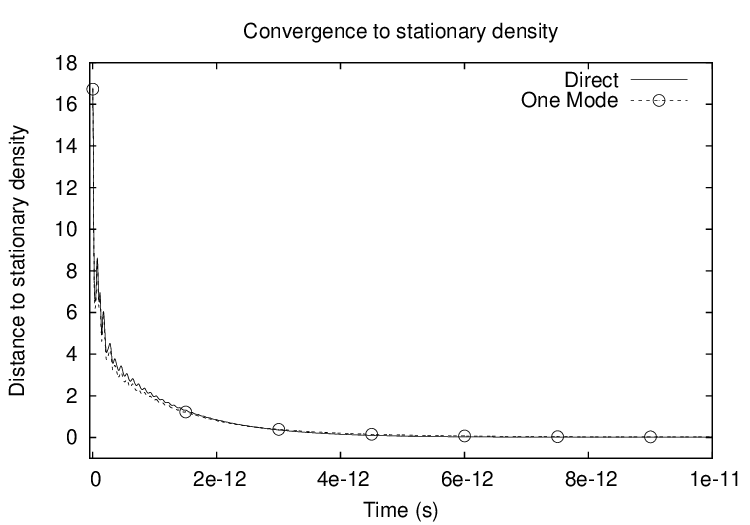}\\
\includegraphics[width=0.49\linewidth]{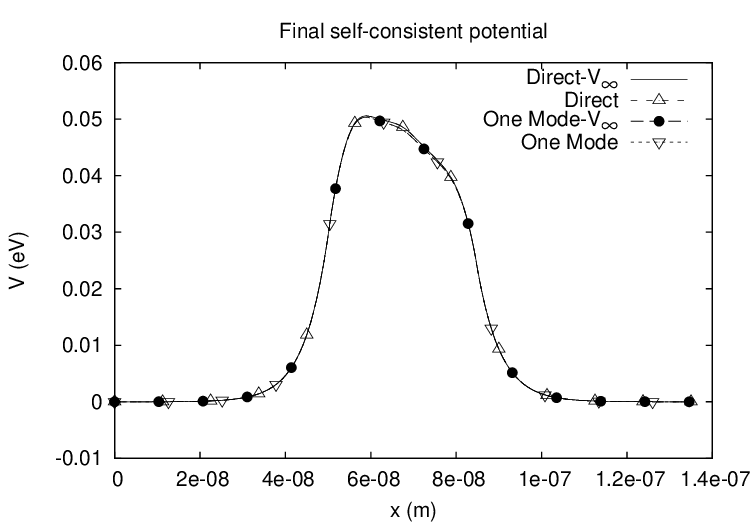}&
\includegraphics[width=0.49\linewidth]{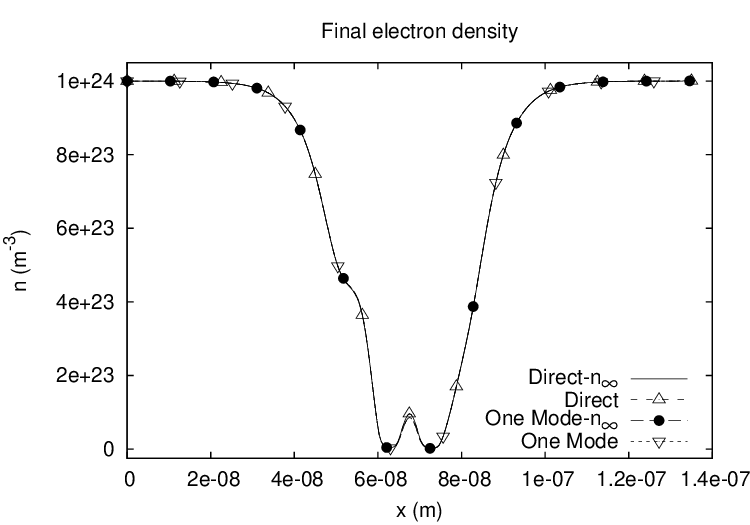}
\end{tabular}
\caption{For $M=500$, $P=3900$, $P'=300$, $\Delta t=2\times10^{-15}\,s$, $l_{tot}=5\times10^{3}$ and, for the methods considered, are depicted from left to right: Evolution of the resonant energy, Evolution of the potential relative distance \eqref{eq_distv} to the stationary solution, Evolution of the average current density defined by \eqref{eq_defcour}, Evolution of the density relative distance \eqref{eq_ecL2} to the stationary solution, Self-consistent potential at final iteration $l_{tot}$, and Electron density at final iteration $l_{tot}$.}\label{graph_t_001}
\end{center}
\end{figure}

We consider now Figure \ref{graph_t_001} and Table \ref{tab2meth_inst} where is given a numerical comparison of the methods presented in section \ref{sec_instat}: the Direct Resolution and the One Mode Approximation. The Direct Resolution plays the role of the reference method.\\
For the Direct Resolution, resp. One Mode Approximation, the initial data of the method and the stationary potential $V_{\infty}$ and density $n_{\infty}$ corresponding to the final bias $B_{\infty}$ are obtained as in section \ref{sec_parInst}, resp. \ref{sec_parInstr}, where $P=3900$ and $P'=300$.\\ 
For the two methods, the boundary conditions for the time-dependent Schr\"odinger equation are given by \eqref{dtbcNHg_cut}\eqref{dtbcNHd_cut} with $M=500$, the corresponding angular frequencies by \eqref{eqwgd}\eqref{eqwdd} and the extrapolation of the intermediary potential by \eqref{eqExtV}. The frequency mesh is given by the integers $P=3900$ and $P'=300$, the total number of iterations is equal to $l_{tot}=5\times10^{3}$ and the time step is equal to $\Delta t = 2\times10^{-15}\,s$. The last verifies the condition \eqref{eqstabt_toutk}. For the One Mode Approximation the angular frequency $\omega^{\infty}_p$ is given by \eqref{eqwrd}.\\
In the graphic of Figure \ref{graph_t_001} entitled {\it Convergence to stationary density}, we depicted the time evolution of the $l^2$-norm relative distance to the density $n_{\infty}$:
\begin{equation}\label{eq_ecL2}
100\frac{|n^l-n_{\infty}|_2}{|n_{\infty}-n_D|_2}\,,
\end{equation}
where $n^{l}$ is the electron density corresponding to the iteration number $l$ of Algorithm \ref{AlgTra}. Since the relevant variation of the density occurs inside the device, we considered in \eqref{eq_ecL2} the relative distance from $(n^l-n_D)$ to $(n_{\infty}-n_D)$. In the graphic entitled {\it Evolution of the resonant energy}, we depicted the time evolution of the real part of the first resonance $z^l$ computed using the method presented in section \ref{sec_calres} where the potential $Q$ is equal to $U^l+V^l$. In the last two graphics of Figure \ref{graph_t_001}, we represented the self-consistent potential and the probability density at final time step provided by the Direct Resolution, curve titled Direct, and the One Mode Approximation, curve titled One Mode. They are compared with the stationary Direct Resolution potential $V_{\infty}$, curve titled Direct-$V_{\infty}$, and density $n_{\infty}$, curve titled Direct-$n_{\infty}$. The same was done for the One Mode Approximation: curves titled One Mode-$V_{\infty}$ and One Mode-$n_{\infty}$. In Table \ref{tab2meth_inst}, the column {\it Final dist} is the value of the quantity \eqref{eq_distv} for $l=l_{tot}$.\\
The two methods behave the same qualitatively: the potential (resp. density) converges, up to a small error, to $V_{\infty}$ (resp. $n_{\infty}$). At time $t=0^+$, the current is close to the current corresponding to $B_I$. This is not true for the energy due to the discontinuity of the bias at $t=0$. After fast transient oscillations, the current stabilizes in a short time to a value relatively close to the current corresponding to $B_{\infty}$ whereas the energy increases slowly to such a value after slightly longer transient oscillations. The comparison is also favorable quantitatively. Indeed, at each time step, the currents provided by the two methods are relatively close of each other which, from formula \eqref{eq_defcour}, signifies that the One Mode Approximation wave functions are computed accurately in comparison with the Direct Resolution. Similarly, at each time step, the small relative difference between the energies provided by the two methods validates the One Mode Approximation self-consistent potential in comparison with the Direct Resolution. We conclude that for the parameters chosen in this section, the One Mode Approximation provides a solution accurate enough at each time step with a reduction of the CPU time by a factor close to $6$ compared to the Direct Resolution, see Table \ref{tab2meth_inst}.\\
\begin{table}
\begin{center}
\begin{tabular}{|c|c|c|c|}
      \hline 
      & CPU($s$) & Final dist ($\%$) \\
      \hline 
      Direct & $952.35$ & $4.4593\times10^{-3}$ \\ 
      \hline
      One Mode & $159.30$ & $1.0857\times10^{-2}$ \\ 
      \hline
    \end{tabular}
\caption{Comparison of the Direct Resolution and One Mode Approximation for the resolution of the time-dependent Schr\"odinger-Poisson system \eqref{schInst}\eqref{eqPot}\eqref{densIns} for $M=500$, $P=3900$, $P'=300$, $\Delta t=2\times10^{-15}\,s$ and $l_{tot}=5\times10^{3}$.}\label{tab2meth_inst}
\end{center}
\end{table}
The above comparison was not performed with the parameters of the reference method which give the best accuracy/CPU time ratio. However, when the parameters of the Direct Resolution are adjusted, as explained in section \ref{sec_parInst}, to reduce the CPU time within a reasonable range of accuracy, the One Mode Approximation provides a solution in the same range of accuracy with a lower CPU time.
\begin{remark}
As it is proposed in \cite{BoFaNi}, the average current density at time $t$ is defined by:
\begin{equation}\label{eq_defcour}
\mathcal{J}(t) = \frac{q\hbar}{m}\int_{\R}g(k)\textrm{Im}\left[\frac{1}{L}\int_0^L\partial_x\Psi_k(t,x)\overline{\Psi_k(t,x)}dx\right]dk\,.
\end{equation}
In \eqref{eq_defcour}, the space derivative is computed using the finite difference method. For the Direct Resolution, the integral with respect to the frequency $k$ is computed using the trapezoidal rule where the frequency mesh is the mesh used to compute the density. For the One Mode Approximation, the frequency integral is given by a double scale trapezoidal rule of the form \eqref{eq_dbscint}. As we did for the density, it is obtained by injecting in \eqref{eq_defcour} the decomposition \eqref{dec_psiIns} of the wave function and neglecting the cross term between the resonant and the non resonant part.
\end{remark} 

\subsubsection{Time evolution of the resonance peaks}\label{sec2pics}

\begin{figure}
\begin{center}
\begin{tabular}{cc}
\includegraphics[width=0.5\linewidth]{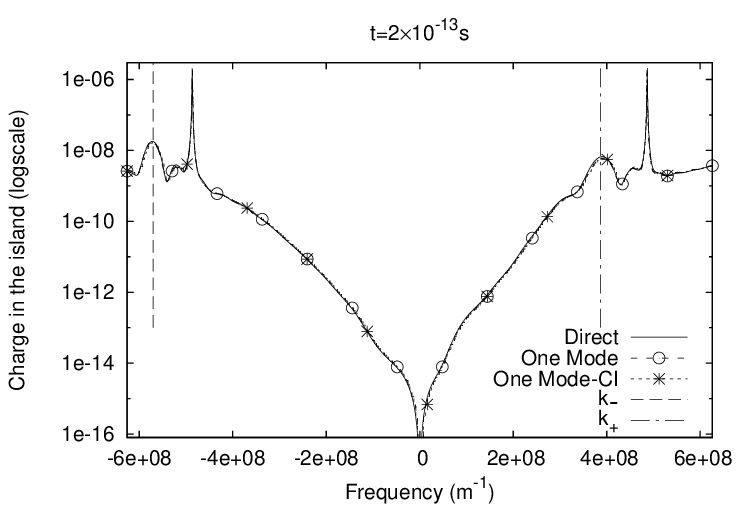}&\includegraphics[width=0.5\linewidth]{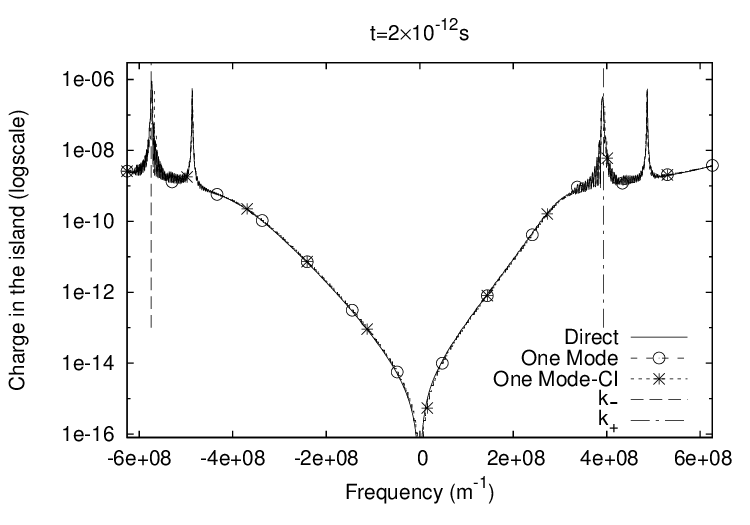}\\
\includegraphics[width=0.5\linewidth]{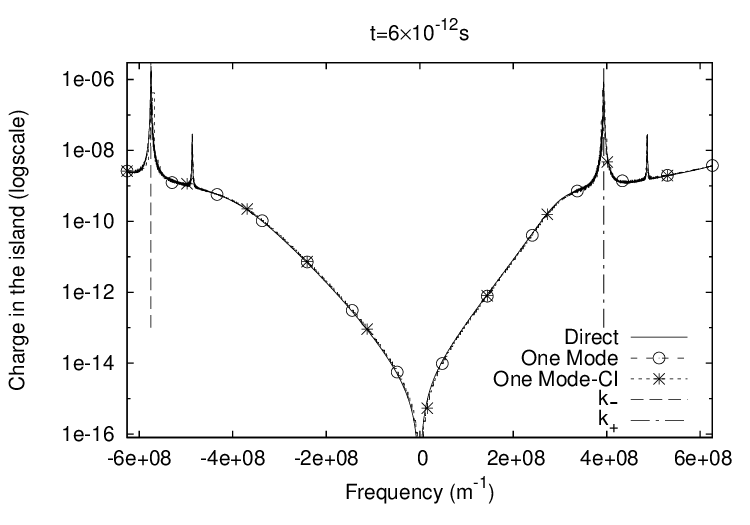}&\includegraphics[width=0.5\linewidth]{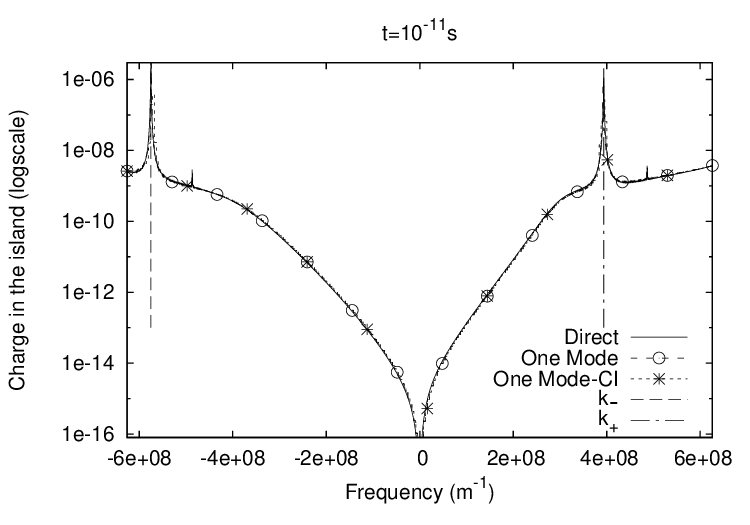}\\
\includegraphics[width=0.5\linewidth]{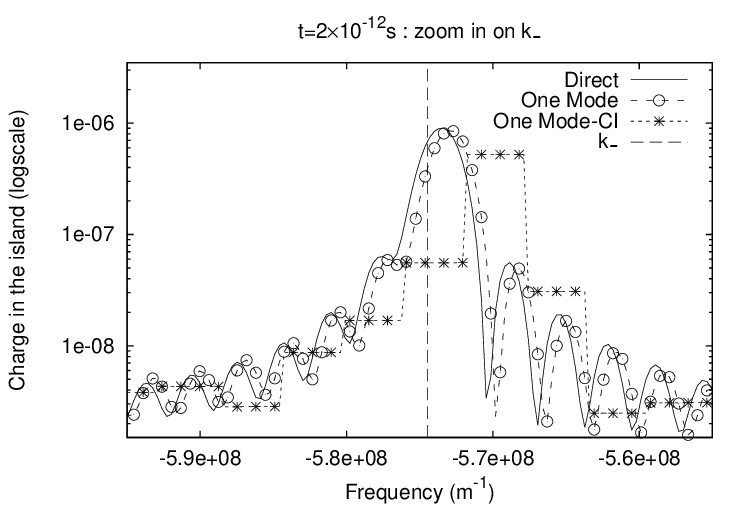}&\includegraphics[width=0.5\linewidth]{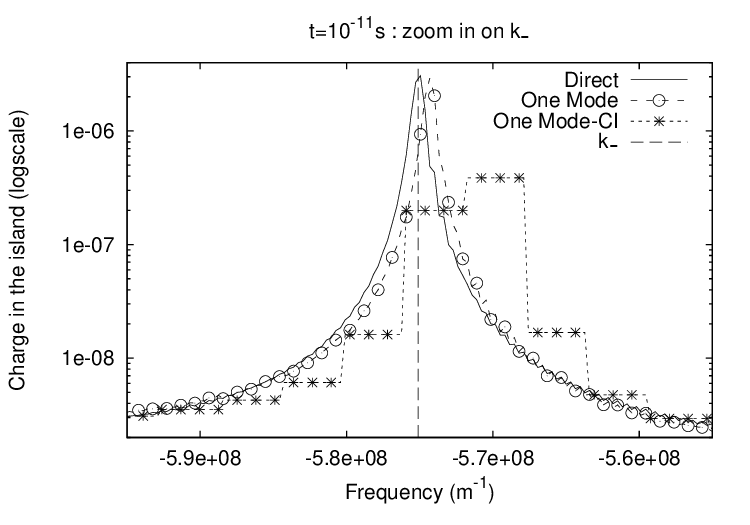}
\end{tabular}
\caption{For $P=3900$, $P'=300$, $\Delta t=10^{-15}\,s$ and at different values of the time $t$: Number \eqref{eq_Nk} of electrons in the interval $[a_2,b_2]$ for one wave function  with respect to the frequency $k$ for the Direct Resolution, Number \eqref{eq_Nkr} of electrons in the interval $[a_2,b_2]$ for one wave function with respect to the frequency $k$ for the One Mode Approximation, and Resonant frequencies given by \eqref{eq_kR} (a vertical line is represented at these frequencies). In the two last graphics, is depicted a zoom in on $k^l_{-}$ of the preceding curves at time $t=2\times10^{-12}\,s$ and $t=10^{-11}\,s$.}\label{evpic_Vet}
\end{center}
\end{figure}

For biases of the form \eqref{B_marche}, since the change of potential at time $t=0$ is abrupt, the adiabaticity hypothesis of \cite{PrSj} is not satisfied and one expects two peaks corresponding to both the resonant energy at time $t=0^-$ and the resonant energy at time $t$. This is what we verify numerically in the present section by using the Direct Resolution and the One Mode Approximation.\\
The parameters and the initialization of the two methods are the same than in section \ref{sec_compt}. Since we are not interested in the reduction of the CPU time, we took the exact boundary conditions \eqref{dtbcNHg}\eqref{dtbcNHd} for the time-dependent Schr\"odinger equation, instead of the simplified ones \eqref{dtbcNHg_cut}\eqref{dtbcNHd_cut}, and the time step equal to $\Delta t =10^{-15}\,s$, instead of $\Delta t = 2\times10^{-15}\,s$.\\ 
In Figure \ref{evpic_Vet}, for the curve with title Direct, the number of electrons in the interval $[a_2,b_2]$ carried by one wave function is given by: 
\begin{equation}\label{eq_Nk}
\displaystyle N^l_{di,p}=\sum_{j\in \tilde{w}}|\Psi^l_{p,j}|^2\Delta x\,,
\end{equation}
where $\Psi^l_{p,j}$ is the Direct Resolution wave function defined in section \ref{sec_algInst} and $\tilde{w}=\{j\,|\,a_2\leq x_j<b_2\}$. For the curve with title One Mode, the cross term between the resonant and the non resonant part of the wave function is neglected and equation \eqref{eq_Nk} becomes:
\begin{equation}\label{eq_Nkr}
N^l_{om,p}=\sum_{j\in \tilde{w}}\left|\Psi^{nr,l}_{p',j}\right|^2\Delta x + \left|\lambda^l_p\right|^2\sum_{j\in \tilde{w}}\left|u^{l-\frac{1}{2}}_j\right|^2\Delta x\,,
\end{equation}
where $\Psi^{nr,l}_{p',j}$, $\lambda^l_p$ and $u^{l-\frac{1}{2}}_j$ are the components, defined in section \ref{sec_S2I}, of the One Mode Approximation wave function. Like in section \ref{sec_S2I}, $0\leq p'\leq P'$ is the integer such that $p=\nu p'+r$ for some $0\leq r\leq \nu-1$, in other terms, we made in \eqref{eq_Nkr} a constant interpolation of the non resonant wave function for $p\neq \nu p'$. Consider the numbers $E^l\in\R$ and $\Gamma^l>0$ such that
\[
z^l=E^l-i\Gamma^l\,,
\]
where $z^l$ is the Direct Resolution resonance, at time $t^l$, defined in section \ref{sec_compt}. Then, the resonant frequencies
\begin{equation}\label{eq_kR}
\displaystyle k^l_{-} = -\sqrt{\frac{2m}{\hbar^2}\left(E^l+B_{\infty}\right)}\,, \quad\quad\quad k^l_{+} = \sqrt{\frac{2m}{\hbar^2}E^l}
\end{equation}
are the real numbers related to the resonant energy $E^l$ according to the relation \eqref{eqDispt}.\\ 
We can now discuss the results in Figure \ref{evpic_Vet}. For the Direct Resolution, we observe that for small times $t^l$, the localization of the resonance peak is not given by the resonant energy $E^l$ but by the resonant energy at time $t=0^-$: $E_I$. Indeed, for such times, the number of electrons $N^l_{di,p}$ defined in \eqref{eq_Nk} has only one peak which is located at the frequencies such that $E_k=E_I$ where $E_k$ and $E_I$ are given by \eqref{eqDisp} and \eqref{eqZI} respectively (the numerical value of $E_I$ is computed as in section \ref{sec_compres} where $B_I=0\,eV$ and $P=3900$). When time goes, the first peak persists but decays slowly while a second peak appears and then grows at the frequencies $k^l_-$ and $k^l_+$ defined in \eqref{eq_kR}. Using the definition of $k^l_-$ and $k^l_+$, the last statement means that the second peak lives at the frequencies such that $\varepsilon_k^\infty=E^l$. At time $t^l=10^{-11}\,s$, the initial peak has almost vanished and the second peak almost reached its maximal amplitude.\\
For the One Mode Approximation, the first peak is provided by the initial condition and, as stated in section \ref{sec_S2I}, the accuracy of the second peak is reached only if the interpolation of the non resonant wave function in \eqref{trapTh} is well chosen. In particular, the constant interpolation computed as explained in section \ref{sec_parInstr}, curve titled One Mode-CI, provides a piecewise constant function close to the resonant frequencies $k^l_-$ and $k^l_+$ and, therefore, does not describe correctly the second peak (this is clear in the two last graphics of Figure \ref{evpic_Vet}). Whereas in the case of the oscillatory interpolation, which corresponds to the iteration \eqref{eqLaint}, the coefficient $\lambda^l_p$ is such that, for large times, the number of electrons $N^l_{om,p}$ defined in \eqref{eq_Nkr} has a sharp peak at the frequencies $k^l_-$ and $k^l_+$  with amplitude high enough and the One Mode Approximation is accurate compared to the Direct Resolution. 

\begin{acknowledgements}\bigskip
The authors acknowledge support from the project QUATRAIN (BLAN07-2 212988) funded by the French Agence Nationale de la Recherche and from the Marie Curie
Project DEASE: MEST-CT-2005-021122 funded by the European Union. The second author acknowledges Prof. Norbert Mauser for his hospitality in WPI Vienna where a part of the present work was achieved.
\end{acknowledgements}

\appendix{}

\section{Resolution of the stationary Schr\"odinger equation using finite difference}\label{app_schSt}

Discrete transparent boundary conditions are proposed in \cite{Ar} to solve equation \eqref{VpGe} on the domain $[0,L]$ using finite difference. In order to be able to justify the choice of the transparent boundary conditions in section \ref{app_RDF}, we recall how the method in \cite{Ar} is deduced under the assumption $k>0$ and we give an extension of this method to the case $k<0$.\\
For a given potential $Q_I$ such that:
\begin{equation}\label{eqHQ}
\displaystyle Q_I(x)=0\,,\, x\leq 0 \quad \textrm{ and } \quad  Q_I(x)=Q_{I,L}\,, \, x\geq L
\end{equation}
and for a given wave vector $k\in\R$, we consider the stationary Schr\"odinger equation
\begin{equation}\label{schstGen}
-\frac{\hbar^2}{2m}\partial_x^2\Phi+Q_I\Phi=E_k\Phi\,, \quad x \in \R
\end{equation} 
with the $k$-dependent scattering conditions \eqref{eqSCkp}, \eqref{eqSCkn}.\\
Using the space discretization defined in section \ref{sec_algStat} and noting $\Phi_j$ the approximation of the wave function $\Phi(x_j)$ and $Q_{I,j}$ the approximation of the potential $Q_I(x_j)$, the finite difference method for equation \eqref{schstGen} writes:
\begin{equation}\label{eqdfk}
D_x^2\Phi_j = \frac{2m}{\hbar^2}(Q_{I,j}-E_k)\Phi_j\,, \quad j=1,...,J-1\,, 
\end{equation}
where
\[
D_x^2\Phi_j = \frac{\Phi_{j+1}-2\Phi_j+\Phi_{j-1}}{\Delta x^2}\,.
\]
For $k>0$, the energy is given by $E_k=\frac{\hbar^2k^2}{2m}$. Then, thanks to assumption \eqref{eqHQ}, the solutions to equation \eqref{eqdfk} can be computed explicitly for $j\leq 1$. Indeed, under the condition 
\begin{equation}\label{eqstabx}
\Delta x < \frac{1}{\sqrt{k^2+\frac{2m}{\hbar^2}\vert Q_{I,L}\vert}}\,,
\end{equation}
they correspond to the linear combinations of $\alpha_1^j$, $\alpha_2^j$ where:
\begin{equation}\label{eqal12}
\alpha_{1,2}= 1-\Delta x^2\frac{k^2}{2}\pm i\sqrt{\Delta x^2k^2-\Delta x^4\frac{k^4}{4}}
\end{equation}
are complex numbers verifying $\vert\alpha_j\vert=1$ and $\alpha_1=e^{i\tilde{k}\Delta x}$ for some $\tilde{k}>0$. Then $\alpha_1^j$ is an incoming plane wave, $\alpha_2^j=\alpha_1^{-j}$ is the reflected wave and the discrete analogous to the scattering condition \eqref{eqSCkp} is to select the linear combinations of $\alpha_1^j$, $\alpha_2^j$ of the form:
\begin{equation}\label{eqPhidg}
\Phi_j = \alpha_1^{j} + R \alpha_1^{-j}\,, \quad j \leq 1\,.
\end{equation}
Applying equation \eqref{eqPhidg} to compute $\Phi_0$ and $\Phi_1$, we get the following discrete transparent boundary condition at $x=0$:
\begin{equation}\label{eqBCkpg}
\alpha_1^{-1}\Phi_0-\Phi_1=\alpha_1^{-1}-\alpha_1\,.
\end{equation}
We note that inequality \eqref{eqstabx} implies $\tilde{k}\Delta x<\frac{\pi}{2}$ and the discretization provides enough space points in one wave length of the discrete plane wave.\\
Similarly, for $j\geq J$ it holds $Q_{I,j}=Q_{I,L}$ and, if $k^2-\frac{2m}{\hbar^2}Q_{I,L}>0$, the solutions to equation \eqref{eqdfk} are the linear combinations of $\beta_1^j$, $\beta_2^j$ where:
\begin{equation}\label{eqbe12kp}
\beta_{1,2}= 1-\Delta x^2(\frac{k^2}{2}-\frac{m}{\hbar^2}Q_{I,L}) \pm i \sqrt{\Delta x^2(k^2-\frac{2m}{\hbar^2}Q_{I,L})-\Delta x^4(\frac{k^2}{2}-\frac{m}{\hbar^2}Q_{I,L})^2}\,.
\end{equation}
Inequality \eqref{eqstabx} implies $\vert\beta_j\vert=1$ and $\beta_1=e^{i\tilde{k}\Delta x}$ for some $\tilde{k}>0$. Then $\beta_1^j$ is an outgoing wave and the scattering condition \eqref{eqSCkp} gives:
\begin{equation}\label{eqPhidd}
\Phi_j = T\beta_1^{j}\,, \quad j \geq J-1\,.
\end{equation}
Applying equation \eqref{eqPhidd} to compute $\Phi_{J-1}$ and $\Phi_J$, we get the following discrete transparent boundary condition at $x=L$:
\begin{equation}\label{eqBCkpd}
-\Phi_{J-1}+\beta_1^{-1}\Phi_J=0\,.
\end{equation}
Again, condition \eqref{eqstabx} implies $\tilde{k}\Delta x<\frac{\pi}{2}$. If $k^2-\frac{2m}{\hbar^2}Q_{I,L}\leq 0$, the bounded solutions to \eqref{eqdfk} for $j\geq J-1$ are given by \eqref{eqPhidd}, where $\Phi_j$ is constant in the case of the equality and vanishing in the case of the strict inequality, and \eqref{eqBCkpd} still apply.\\
For $k<0$, the energy is given by $E_k=\frac{\hbar^2k^2}{2m}+Q_{I,L}$. Then, using the assumptions \eqref{eqHQ}, \eqref{eqstabx} and repeating the above calculations where the scattering condition \eqref{eqSCkp} is replaced by \eqref{eqSCkn}, we obtain the following transparent boundary conditions for equation \eqref{eqdfk} at $x=0$ and $x=L$:
\begin{equation}\label{eqBCkng}
\beta_2\Phi_{0}-\Phi_1=0
\end{equation}
\begin{equation}\label{eqBCknd}
-\Phi_{J-1}+\alpha_2\Phi_J=\alpha_2^J(\alpha_2-\alpha_2^{-1})\,,
\end{equation}
where $\alpha_2$ is given by equation \eqref{eqal12} and 
\begin{equation}\label{eqbe12kn}
\beta_{1,2}= 1-\Delta x^2(\frac{k^2}{2}+\frac{m}{\hbar^2}Q_{I,L}) \pm i \sqrt{\Delta x^2(k^2+\frac{2m}{\hbar^2}Q_{I,L})-\Delta x^4(\frac{k^2}{2}+\frac{m}{\hbar^2}Q_{I,L})^2}\,.
\end{equation} 
For $k>0$, the scheme \eqref{eqdfk} with the boundary conditions \eqref{eqBCkpg}\eqref{eqBCkpd} has a unique solution. By construction, the unique solution $\Phi_j$ to the whole-space scheme \eqref{eqdfk} considered on $j\in\mathbb{Z}$ with the scattering conditions \eqref{eqPhidg}\eqref{eqPhidd} verifies the problem \eqref{eqdfk}\eqref{eqBCkpg}\eqref{eqBCkpd}. It follows that the solution to \eqref{eqdfk}\eqref{eqBCkpg}\eqref{eqBCkpd} corresponds exactly with the restriction of $\Phi_j$ to $0\leq j\leq J$. A similar property holds for $k<0$. Moreover we have the following estimate on the solution:
\begin{lemma}\label{le_estst}
There exists a constant $C>0$ which depends only on $\hbar$, $m$, $L$, $Q_{I,L}$, $\sup_{x\in(0,L)}\vert Q_I(x)\vert$, $k$ such that for all $\Delta x$ verifying \eqref{eqstabx}, the solution to the scheme \eqref{eqdfk} with the boundary conditions \eqref{eqBCkpg}\eqref{eqBCkpd} for $k>0$ and \eqref{eqBCkng}\eqref{eqBCknd} for $k<0$ verifies:
\[
\Vert \Phi \Vert_2 \leq C\,,
\]
where 
\[
\Vert \Phi \Vert_2^2 = \sum_{j=1}^{J-1}\vert\Phi_j\vert^2\Delta x\,.
\]
\end{lemma}

\section{Computation of resonances using finite difference}\label{app_RDF}

For a given potential $Q$ verifying \eqref{eqHQN}, we give in this section a discrete version of the following problem: find $(u,z)$ solution to \eqref{pbMres} such that $u$ is purely outgoing outside the interval $[0,L]$.\\
Using the decay at infinity of the function $g$ in \eqref{FermiDirac}, it follows that only the resonances with real part smaller than $\frac{\hbar^2\kappa^2}{2m}$ are important for the computation of the density, where $\kappa$ is defined in section \ref{sec_algStat}. Then, if we suppose that the resonances have a positive real part, which is not a restriction in the applications, the problem of the computation of resonances can be restricted to the strip
\[
S = \left\{ z\in\mathbb{C}\,\vert\, 0<\textrm{Re}\,(z)<\frac{\hbar^2\kappa^2}{2m}\right\}\,.
\]
In order to work with analytic functions on $S$, the value of the square root function $\sqrt{z}$ will be the principal value.\\
Proceeding as in section \ref{app_schSt}, we write first the finite difference method for equation \eqref{pbMres}: 
\begin{equation}\label{eqdfr}
D_x^2u_j = \frac{2m}{\hbar^2}(Q_j-z)u_j\,, \quad j=1,...,J-1\,, 
\end{equation}
where $u_j$ denotes the approximation of the resonant mode $u(x_j)$ and $Q_j$ the approximation of the potential $Q(x_j)$. Thanks to assumption \eqref{eqHQN}, the solutions to equation \eqref{eqdfr} can be computed explicitly for $j\leq 1$. Indeed, under the condition
\begin{equation}\label{eqcdxr}
\Delta x < \frac{1}{\sqrt{\kappa^2+\frac{2m}{\hbar^2}\vert Q_L\vert}}\,,
\end{equation}
they correspond to the linear combinations of $\alpha_1(z)^j$, $\alpha_2(z)^j$ where:
\begin{equation}\label{eqal12r}
\alpha_{1,2}(z)= 1-\Delta x^2\frac{m}{\hbar^2}z\pm i\sqrt{2\Delta x^2\frac{m}{\hbar^2}z-\Delta x^4\frac{m^2}{\hbar^4}z^2}\,.
\end{equation}
Utilizing the forthcoming equation \eqref{eq_icoeff}, it holds $\alpha_2(z)=\rho e^{ik\Delta x}$ with $\rho>0$ and $k<0$. Since we are looking for solutions which are outgoing outside $[0,L]$, we impose:
\begin{equation}\label{equdgr}
u_j = R\alpha_2(z)^{j}\,, \quad j \leq 1\,.
\end{equation}
Applying equation \eqref{equdgr} to compute $u_{0}$ and $u_1$, we get the following discrete transparent boundary condition at $x=0$:
\begin{equation}\label{eqBCgr}
\alpha_2(z)u_{0}-u_1=0\,.
\end{equation}
Similarly, for $j\geq J$ it holds $Q_j=Q_L$ and the solutions to equation \eqref{eqdfr} are the linear combinations of $\beta_1(z)^j$, $\beta_2(z)^j$ where:
\begin{equation}\label{eqbe12r}
\beta_{1,2}(z)= 1-\Delta x^2\frac{m}{\hbar^2}(z-Q_L)\pm i\sqrt{2\Delta x^2\frac{m}{\hbar^2}(z-Q_L)-\Delta x^4\frac{m^2}{\hbar^4}(z-Q_L)^2}\,.
\end{equation}
Using equation \eqref{eq_icoeff}, the outgoing wave is identified with the plus sign in \eqref{eqbe12r}, therefore we impose:
\begin{equation}\label{equddr}
u_j = T\beta_1(z)^{j}\,, \quad j \geq J-1\,.
\end{equation}
Applying equation \eqref{equddr} to compute $u_{J-1}$ and $u_J$, we get the following discrete transparent boundary condition at $x=L$:
\[
-u_{J-1}+\beta_1(z)^{-1}u_J=0\,.
\]
Since $\beta_1(z)\beta_2(z)=1$, the previous equation writes:
\begin{equation}\label{eqBCdr}
-u_{J-1}+\beta_2(z)u_J=0\,.
\end{equation}
In order to work with quantities of order $1$, the condition $u^Hu=1$ is imposed instead of the norm condition in \eqref{pbMres} and the resonant mode have to be divided by $\Delta x^{\frac{1}{2}}$ to verify $\sum_{j=0}^{J}\vert u_j\vert^2\Delta x=1$. Then, it follows from equations \eqref{eqdfr}, \eqref{eqBCgr} and \eqref{eqBCdr} that the sought discrete problem is
\begin{equation}\label{eqpbresd}
\displaystyle \left\{\begin{array}{l}M(z)u=0\\[1.65mm]
u^Hu=1
\end{array}\right.\,,
\end{equation}
where
\begin{equation}\label{eqMz}
M(z)=\left(\begin{array}{ccccc} 
\alpha_2(z)&-1& & & \\
-1&2+\frac{2m\Delta x^2}{\hbar^2}(Q_1-z)& -1& & 0\\
& \ddots & \ddots & \ddots & \\
0 &  & -1 &2+\frac{2m\Delta x^2}{\hbar^2}(Q_{J-1}-z)& -1\\
& & & -1& \beta_2(z)
\end{array}\right)\,.
\end{equation} 
We note that the holomorphy of $M(z)$ required in section \ref{sec_calres} is verified in $S$.\\
We give now an interpretation of the solutions to the problem \eqref{eqpbresd}. In the strip $S$, resonances can be defined at the discrete level as the complex numbers $z$ such that the whole-space scheme \eqref{eqdfr} considered on $j\in\mathbb{Z}$ has purely outgoing solutions, i.e. solutions verifying \eqref{equdgr} and \eqref{equddr}. Such a solution is a called a resonant mode. Then the following properties hold (proof left to the reader): 
\begin{itemize}
\item When $z$ is a resonance, the corresponding resonant mode is uniquely determined up to a phase factor.
\item A complex number $z$ is a resonance if and only if $z$ is such that there exists a vector $u$ verifying \eqref{eqpbresd}.
\item If $z$ is a resonance, a corresponding $u$ solution to \eqref{eqpbresd} is equal, up to a phase factor, to the restriction of the resonant mode to $0\leq j\leq J$.
\end{itemize}    
\begin{proposition}\label{Pr_res}
If $\Delta x$ verifies the condition \eqref{eqcdxr}, then for all $z$ in $S$ the coefficients defined in \eqref{eqal12r} and \eqref{eqbe12r} verify:
\begin{equation}\label{eq_icoeff}
\textrm{Im}\,(\alpha_1(z))>0\,, \quad \textrm{Im}\,(\alpha_2(z))<0 \quad \textrm{and} \quad
\textrm{Im}\,(\beta_1(z))>0\,, \quad \textrm{Im}\,(\beta_2(z))<0\,.
\end{equation}
If in addition $\textrm{Im}\,(z)<0$, then it holds:
\begin{equation*}
\vert\alpha_1(z)\vert>1\,, \quad \vert\alpha_2(z)\vert<1 \quad \textrm{and} \quad
\vert\beta_1(z)\vert>1\,, \quad \vert\beta_2(z)\vert<1\,.
\end{equation*}
\end{proposition}
In the next proposition, we show that resonances have a negative imaginary part. This, combined with Proposition \ref{Pr_res} and equations \eqref{equdgr}\eqref{equddr}, implies that the corresponding resonant modes tend to the infinity when $j\rightarrow\pm\infty$. It is the discrete transcription of the fact that the space $L^2(\R)$ has to be deformed in order to consider resonances as eigenvalues. 
\begin{proposition}
If $\Delta x$ verifies the condition \eqref{eqcdxr} and if $z\in S$ is such that the problem \eqref{eqpbresd} has a solution $u$, then it holds $\textrm{Im}\,(z)<0$. 
\end{proposition}
\begin{proof}
We multiply equation \eqref{eqdfr} by $\overline{u_j}$ and sum up for $j=1,...,J-1$. Then, utilizing the summation by part rule:
\begin{equation}\label{eq_sumpart}
\sum_{j=1}^{J-1}g_jD_x^-f_j = -\sum_{j=0}^{J-1}f_jD_x^+g_j+\frac{1}{\Delta x}(f_{J-1}g_{J}-f_0g_0)\,,
\end{equation}
where 
\[
D_x^+f_j=\frac{f_{j+1}-f_j}{\Delta x}\,, \quad D_x^-f_j=\frac{f_j-f_{j-1}}{\Delta x}\,, \quad D_x^-D_x^+f_j=D_x^2f_j\,,
\]
we get
\begin{equation}\label{eq_ippst}
\sum_{j=0}^{J-1}\left\vert D_x^+u_j\right\vert^2=-\sum_{j=1}^{J-1}\frac{2m}{\hbar^2}(Q_j-z)\vert u_j\vert^2+\frac{1}{\Delta x}\left(D_x^+u_{J-1}\overline{u_J}-D_x^+u_0\overline{u_0}\right)\,.
\end{equation}
Inserting the boundary conditions \eqref{eqBCgr}\eqref{eqBCdr} in \eqref{eq_ippst} and taking the imaginary part, it follows:
\[
0 = \sum_{j=1}^{J-1}\frac{2m}{\hbar^2}\textrm{Im}\,(z)\vert u_j\vert^2 - \frac{\textrm{Im}\,(\beta_2(z))}{\Delta x^2}\vert u_J\vert^2 - \frac{\textrm{Im}\,(\alpha_2(z))}{\Delta x^2}\vert u_0\vert^2\,.
\]
Using $u^Hu=1$ and \eqref{eq_icoeff}, the previous equation implies $\textrm{Im}\,(z)<0$.
\end{proof}

\section{Resolution of the time-dependent Schr\"odinger equation}\label{app_schInst}

Although in the applications considered here the potential depends nonlinearly on the wave function, this section is written in the case of the linear Schr\"odinger equation. As in \cite{ChJiSu} and \cite{ZiEh}, the extension to the nonlinear case is provided by an adapted choice of the potential at half time step. 
\subsection{The homogeneous case}\label{sec_hdtbc}

In this section, we recall the scheme proposed in \cite{EhAr} to solve on the bounded domain $[0,L]$ the time-dependent Schr\"odinger equation
\begin{equation}\label{schinstGen}
\displaystyle \left\{\begin{array}{ll}
i\hbar\partial_t\Psi = -\frac{\hbar^2}{2m}\partial_x^2\Psi+Q\Psi\,, &
t>0\,, \, x \in \R\\[1.65mm]
\Psi(0,x) = \Phi(x)\,,& 	x \in \R
\end{array}\right.\,,
\end{equation} 
with the following hypothesis:
\begin{itemize} 
\item[$H1.$]The initial condition $\Phi$ is supported in $0<x<L$.
\item[$H2.$]The potential $Q$ verifies: for $t>0$
\[
\displaystyle Q(t,x)=0\,,\, x\leq 0 \quad \textrm{ and } \quad  Q(t,x)=Q_L\,, \, x\geq L\,.
\]
\end{itemize} 
 Considering the time and space discretization defined in section \ref{sec_algInst}, we note $Q^{l+1/2}_j$ the approximation of the potential $Q(t^{l+1/2},x_j)$ and $\Psi^l_j$ the approximation of the solution $\Psi(t^l,x_j)$. Then, equation \eqref{schinstGen}, is solved with the Crank-Nicolson method:
\begin{equation}\label{CranNich}
\displaystyle i\hbar D_t^+\Psi_j^l=-\frac{\hbar^2}{2m}D_x^2\Psi^{l+1/2}_j+Q^{l+1/2}_j\Psi^{l+1/2}_j\,,\quad j=1,...,J-1\,,\quad l\geq 0\,, 
\end{equation}
where
\begin{equation}\label{eq_defdtde}
D_t^+\Psi_j^l = \frac{\Psi_j^{l+1}-\Psi_j^l}{\Delta t} \quad \textrm{ and } \quad \Psi^{l+1/2}_j=\frac{1}{2}\left( \Psi_j^l+\Psi_j^{l+1} \right)\,.
\end{equation}
Equation \eqref{CranNich} comes with the discrete transparent boundary conditions: 
\begin{equation}\label{DTBCG}
\displaystyle \Psi^l_1 - s_0^0 \Psi^l_0 = \sum_{k=1}^{l-1} s^{l-k}_0 \Psi^k_0 - \Psi^{l-1}_1\,, \quad  l\geq 1
\end{equation}
\begin{equation}\label{DTBCD}
\displaystyle \Psi^l_{J-1} - s_J^0 \Psi^l_J = \sum_{k=1}^{l-1} s^{l-k}_J \Psi^k_J - \Psi^{l-1}_{J-1}\,, \quad  l\geq 1\,,
\end{equation}
where we have for $j=0$ and $j=J$:
\begin{equation}\label{coeffDtbc}
\displaystyle s_j^l = \left[1-i\frac{R}{2}+\frac{\sigma_j}{2}\right]\delta_l^0+\left[1+i\frac{R}{2}+\frac{\sigma_j}{2}\right]\delta_l^1+\alpha_j\exp\left(-il\varphi_j\right)\frac{P_l(\mu_j)-P_{l-2}(\mu_j)}{2l-1}
\end{equation}
and
\[
\displaystyle \varphi_j = \arctan\frac{2R(\sigma_j+2)}{R^2-4\sigma_j-\sigma_j^2}, \quad \mu_j = \frac{R^2+4\sigma_j+\sigma_j^2}{\sqrt{(R^2+\sigma_j^2)[R^2+(\sigma_j+4)^2]}}\,,
\]
\[
\displaystyle \sigma_j = \frac{2m\Delta x^2}{\hbar^2}Q_j\,, \quad
\alpha_j=\frac{i}{2}((R^2+\sigma_j^2)[R^2+(\sigma_j+4)^2])^{1/4}\exp\left(i\frac{\varphi_j}{2}\right)\,, \quad R=\frac{4m\Delta x^2}{\hbar\Delta t}\,.
\]
Here $P_l$ denotes the Legendre polynomials with the convention $P_{-1} = P_{-2} = 0$, $\delta_l^j$ the Kronecker symbol related to the integers $j,\,l$ and $Q_0=0$, $Q_J=Q_L$. The coefficients $s_j^l$ are computed thanks to the recursion formula derived in \cite{EhAr} from the standard recursion formula for the Legendre polynomials.\\
We remark that the scheme \eqref{CranNich} with the boundary conditions \eqref{DTBCG}\eqref{DTBCD} has a unique solution. By construction, see \cite{EhAr}, the unique solution $\Psi_j^l$ in $L^2(\mathbb{Z})$ to the whole-space scheme \eqref{CranNich} considered on $j\in\mathbb{Z}$ verifies the problem \eqref{CranNich}\eqref{DTBCG}\eqref{DTBCD}. It follows that the solution to \eqref{CranNich}\eqref{DTBCG}\eqref{DTBCD} corresponds exactly with the restriction of $\Psi_j^l$ to $0\leq j\leq J$. Moreover, we have the following stability result:
\begin{theorem}[\cite{EhAr}]\label{th_sth}
The solution to the discrete Schr\"odinger equation \eqref{CranNich} with the boundary conditions \eqref{DTBCG}\eqref{DTBCD} is uniformly bounded:
\[
\Vert \Psi^l \Vert_2 \leq \Vert \Psi^0 \Vert_2\,, \quad l\geq 1 
\]
and therefore, the scheme is unconditionally stable.
\end{theorem}

\subsection{The non-homogeneous case}\label{sec_nhdtbc}

\subsubsection{The non-homogeneous discrete transparent boundary conditions}

We consider the problem \eqref{schinstGen} where the initial condition $\Phi$ is solution to \eqref{schstGen}. We suppose that the potentials $Q_I$ and $Q$ are such that \eqref{eqHQ} and $H2$ hold and $Q(0,x)=Q_I(x)$. The initial condition does not verify assumption $H1$, therefore \eqref{DTBCG}\eqref{DTBCD} have to be replaced by suitable boundary conditions in the Crank-Nicolson scheme \eqref{CranNich}.\\  
Using equations \eqref{schinstGen} and \eqref{schstGen}, the function
\begin{equation}\label{eqfdiffg}
\displaystyle \varphi = \Psi - \Phi e^{-i\omega_0t}\,,
\end{equation}
where
\begin{equation}\label{eqwg}
\omega_0 = \frac{E_k}{\hbar}\,,
\end{equation}
is solution to:
\begin{equation*}
\displaystyle i\hbar\partial_t\varphi  = [-\frac{\hbar^2}{2m}\partial_x^2+Q]\varphi\,, \quad  x\leq 0
\end{equation*}
and verifies $\varphi(0,x)=0$. Thus, we can write the homogeneous boundary condition \eqref{DTBCG} for $\varphi$ and obtain the following boundary condition at $x=0$ for $\Psi$:
\begin{multline}\label{dtbcNHg}
\displaystyle \Psi^l_1 - s^0_0\Psi^l_0 = \sum_{k=1}^{l-1}s^{l-k}_0\left(\Psi^k_0- \Phi_0e^{-i\omega_0t^k}\right)-\left(\Psi^{l-1}_1-\Phi_1e^{-i\omega_0t^{l-1}}\right) \\
+ \Phi_1e^{-i\omega_0t^l} - s^0_0\Phi_0e^{-i\omega_0t^l}\,,\quad l\geq 1\,,
\end{multline}
where the coefficients $s_j^l$ are given by \eqref{coeffDtbc} and $\Phi_j$ is a solution to the stationary discrete problem \eqref{eqdfk}. We proceed similarly at $x=L$ by setting
\begin{equation}\label{eqfdiffd}
\displaystyle \varphi = \Psi - \Phi e^{-i\omega_Lt}\,,
\end{equation}
where
\begin{equation}\label{eqwd}
\omega_L = \frac{1}{\hbar}(E_k + Q_L - Q_{I,L})\,.
\end{equation}
This leads to the following boundary condition at $x=L$:
\begin{multline}\label{dtbcNHd}
\displaystyle \Psi^l_{J-1} - s^0_J\Psi^l_J = \sum_{k=1}^{l-1}s^{l-k}_J\left(\Psi^k_J-\Phi_Je^{-i\omega_Lt^k}\right) - \left(\Psi^{l-1}_{J-1}-\Phi_{J-1}e^{-i\omega_Lt^{l-1}}\right)\\
+ \Phi_{J-1}e^{-i\omega_Lt^l} - s^{0}_J\Phi_Je^{-i\omega_Lt^l}\,, \quad l\geq 1\,.
\end{multline}

\subsubsection{Discrete angular frequencies, stability}\label{sec-discangf}

Since the homogeneous transparent boundary conditions \eqref{DTBCG}\eqref{DTBCD} are derived from the fully discrete scheme \eqref{CranNich} where $j\in\mathbb{Z}$, the accuracy of the non-homogeneous boundary conditions \eqref{dtbcNHg}\eqref{dtbcNHd} is improved under the condition that the functions in \eqref{eqfdiffg} and \eqref{eqfdiffd} are solutions to the discrete equation \eqref{CranNich} for $j\leq 0$ and, respectively, $j\geq J$. As shown in \cite{Ar}, this leads to the non-homogeneous boundary conditions \eqref{dtbcNHg}\eqref{dtbcNHd}, where the angular frequencies given in \eqref{eqwg} and \eqref{eqwd} are replaced by the discrete ones below:
\begin{equation}\label{eqwgd}
\omega_0 = \frac{2}{\Delta t}\arctan\left(\frac{\Delta t E_k}{2\hbar}\right)\,,
\end{equation}
\begin{equation}\label{eqwdd}
\omega_L = \frac{2}{\Delta t}\arctan\left(\frac{\Delta t}{2\hbar}(E_k + Q_L - Q_{I,L})\right)\,.
\end{equation}
We remark that, for $k>0$, the solution to the initial finite difference scheme \eqref{eqdfk} considered on $j\in\mathbb{Z}$ with the scattering conditions \eqref{eqPhidg}\eqref{eqPhidd} is in $L^{\infty}(\mathbb{Z})$ but not in $L^2(\mathbb{Z})$. The same property is true for $k<0$. Therefore, for such an initial condition, the suitable state space for the solutions to the scheme \eqref{CranNich} considered on $j\in\mathbb{Z}$ is $L^2(\mathbb{Z})+L^{\infty}(\mathbb{Z})$. In particular, it is easy to show that this scheme has a unique solution $\Psi_j^l$ such that $\left(\Psi_j^l\right)_{j\in\mathbb{Z}}\in L^2(\mathbb{Z})+L^{\infty}(\mathbb{Z})$ for $l \geq 0$. By construction, $\Psi_j^l$ verifies the problem \eqref{CranNich} with the boundary conditions \eqref{dtbcNHg}\eqref{dtbcNHd} where the angular frequencies are given by \eqref{eqwgd}\eqref{eqwdd}. It follows that the unique solution to \eqref{CranNich}\eqref{dtbcNHg}\eqref{dtbcNHd} corresponds exactly with the restriction of $\Psi_j^l$ to $0\leq j\leq J$. Moreover, we have the following stability result:
\begin{theorem}\label{th_stnh}
Let $\Phi$ be the solution to the scheme \eqref{eqdfk} with the boundary conditions \eqref{eqBCkpg}\eqref{eqBCkpd} for $k>0$ and \eqref{eqBCkng}\eqref{eqBCknd} for $k<0$. Let $\Psi$ be the solution to the discrete Schr\"odinger equation \eqref{CranNich} with the non-homogeneous boundary conditions \eqref{dtbcNHg}\eqref{dtbcNHd} and the initial condition $\Psi^0=\Phi$. We suppose in addition that the angular frequencies $\omega_0$ and $\omega_L$ in \eqref{dtbcNHg}\eqref{dtbcNHd} are given by \eqref{eqwgd}\eqref{eqwdd}. Then, there exists a constant $C>0$ which depends only on $\hbar$, $m$, $L$, $Q_{I,L}$, $\sup_{x\in(0,L)}\vert Q_I(x)\vert$, $Q_L$, $k$ such that for all $\Delta x$ verifying \eqref{eqstabx} it holds for all $N\geq 0$:
\begin{equation}\label{eq_estth}
\Vert \Psi^l \Vert_2 \leq C\left[\left(1+\max_{0\leq n\leq N}\Vert Q^{n+1/2} \Vert_{\infty}\right)T+1\right]\,, \quad 0\leq l\leq N\,, 
\end{equation}
where $T=N\Delta t$ and $\Vert Q^{n+1/2} \Vert_{\infty} = \max_{1\leq j \leq J-1}\vert Q^{n+1/2}_j \vert$.
\end{theorem}
Even if no size limit on the time step is deduced from Theorem \ref{th_stnh}, the time oscillating terms $e^{-i\omega_0t}$ and $e^{-i\omega_Lt}$ in \eqref{dtbcNHg} and, respectively, \eqref{dtbcNHd} require the following condition for a reasonable resolution:
\begin{equation}\label{eqstabt}
\vert \omega_j\vert\Delta t < \frac{\pi}{2}\,, \quad j=0,L\,.
\end{equation}
We can now give the proof of Theorem \ref{th_stnh}.\\
\begin{proof}
Without any practical restriction, we can suppose that $J\geq 3$. Let us consider a smooth function $\chi$ such that 
\[
0\leq \chi(x)\leq 1\,,\;x\in\R\,,\quad \chi(x)=1\,,\;x\leq\frac{L}{3}\,,\quad \chi(x)=0\,,\;x\geq\frac{2L}{3}\,.
\]
Using the vector $\chi_j=\chi(x_j)$, we define for $l\geq 0$ and $0\leq j\leq J$:
\[
\theta_j^l = \chi_j\Phi_je^{-i\omega_0t^l} + (1-\chi_j)\Phi_je^{-i\omega_Lt^l} 
\]  
and
\[
\varphi_j^l =\Psi_j^l - \theta_j^l\,.
\]
Using the scheme \eqref{CranNich} for $\Psi$, we obtain the discrete Schr\"odinger equation below for the function $\varphi$:
\begin{equation}\label{eq_CNST}
i\hbar D_t^+\varphi_j^l=-\frac{\hbar^2}{2m}D_x^2\varphi^{l+1/2}_j+Q^{l+1/2}_j\varphi^{l+1/2}_j + F_j^l\,,\quad j=1,...,J-1\,,\quad l\geq 0
\end{equation}
with source term
\[
 F_j^l = -i\hbar D_t^+\theta_j^l-\frac{\hbar^2}{2m}D_x^2\theta^{l+1/2}_j+Q^{l+1/2}_j\theta^{l+1/2}_j\,,
\]
where $\varphi^{l+1/2}_j$, $\theta^{l+1/2}_j$ are defined as in \eqref{eq_defdtde}. Multiplying equation \eqref{eq_CNST} by $\overline{\varphi}_j^{l+1/2}$ and summing up for $j=1,...,J-1$, it follows:
\begin{multline*}
i\hbar \sum_{j=1}^{J-1}D_t^+\varphi_j^l\overline{\varphi}_j^{l+1/2}= \frac{\hbar^2}{2m}\sum_{j=0}^{J-1}\left\vert D_x^+\varphi_j^{l+1/2}\right\vert^2+\sum_{j=1}^{J-1}Q_j^{l+1/2}\left\vert \varphi_j^{l+1/2}\right\vert^2 \\-\frac{\hbar^2}{2m\Delta x}\left(D_x^+\varphi_{J-1}^{l+1/2}\overline{\varphi}_J^{l+1/2}-D_x^+\varphi_0^{l+1/2}\overline{\varphi}_0^{l+1/2}\right) + \sum_{j=1}^{J-1}F_j^l\overline{\varphi}_j^{l+1/2}\,,
\end{multline*}
where we used the summation by part rule \eqref{eq_sumpart}. Taking the imaginary part, we get:
\[
D_t^+\Vert\varphi^l\Vert_2^2 = \frac{\hbar}{m}\left[-\textrm{Im}\,\left(D_x^-\varphi_J^{l+1/2}\overline{\varphi}_J^{l+1/2}\right)+\textrm{Im}\,\left(D_x^+\varphi_0^{l+1/2}\overline{\varphi}_0^{l+1/2}\right)\right] + \frac{2\Delta x}{\hbar}\sum_{j=1}^{J-1}\textrm{Im}\,\left(F_j^l\overline{\varphi}_j^{l+1/2}\right)\,.
\]
The initial condition $\varphi^0=0$ obviously holds, therefore, we have by summing up with respect to $l$ in the previous equation:
\begin{multline}\label{eq_summa}
\Vert\varphi^{N+1}\Vert_2^2 = \frac{\hbar\Delta t}{m}\left[-\textrm{Im}\,\left(\sum_{l=0}^{N}D_x^-\varphi_J^{l+1/2}\overline{\varphi}_J^{l+1/2}\right)+\textrm{Im}\,\left(\sum_{l=0}^{N}D_x^+\varphi_0^{l+1/2}\overline{\varphi}_0^{l+1/2}\right)\right]\\ + \frac{2\Delta x\Delta t}{\hbar}\sum_{l=0}^{N}\sum_{j=1}^{J-1}\textrm{Im}\,\left(F_j^l\overline{\varphi}_j^{l+1/2}\right)
\end{multline}
for all $N\geq 0$. Since the non-homogeneous boundary conditions \eqref{dtbcNHg}\eqref{dtbcNHd} hold for $\Psi$, it follows that the homogeneous ones \eqref{DTBCG}\eqref{DTBCD} hold for $\varphi$. Therefore, the boundary terms in \eqref{eq_summa} behave like in Theorem \ref{th_sth}. In particular, it is shown in \cite{EhAr} that the boundary conditions \eqref{DTBCG}\eqref{DTBCD} imply:
\[
-\textrm{Im}\,\left(\sum_{l=0}^{N}D_x^-\varphi_J^{l+1/2}\overline{\varphi}_J^{l+1/2}\right)\, \leq \, 0\,, \quad \textrm{Im}\,\left(\sum_{l=0}^{N}D_x^+\varphi_0^{l+1/2}\overline{\varphi}_0^{l+1/2}\right) \, \leq \, 0\,.
\]
Then, it follows from \eqref{eq_summa} that:
\begin{equation}\label{eq_estphiNL}
\Vert\varphi^{N+1}\Vert_2^2 \leq \frac{2\Delta t}{\hbar}\sum_{l=0}^{N}\left\Vert F^l\right\Vert_2\left\Vert\varphi^{l+1/2}\right\Vert_2\,.
\end{equation}
Using the equations \eqref{eqwgd}\eqref{eqwdd} and the equation \eqref{eqdfk}, with the boundary conditions \eqref{eqBCkpg}\eqref{eqBCkpd} for $k>0$ and \eqref{eqBCkng}\eqref{eqBCknd} for $k<0$, and using the result of Lemma \ref{le_estst}, the following estimate holds for the source term:
\begin{equation}\label{eq_estF}
\left\Vert F^l\right\Vert_2 \leq C\left(1+\max_{0\leq n\leq N}\Vert Q^{n+1/2} \Vert_{\infty}\right)\,, \quad 0\leq l\leq N\,,
\end{equation}
where $C$ denotes a constant which depends only on $\hbar$, $m$, $L$, $Q_{I,L}$, $\sup_{x\in(0,L)}\vert Q_I(x)\vert$, $Q_L$, $k$. Then, after an easy computation, the inequalities \eqref{eq_estphiNL} and \eqref{eq_estF} lead to:
\begin{equation}\label{eq_estphi}
\left\Vert \varphi^l\right\Vert_2 \leq C\left(1+\max_{0\leq n\leq N}\Vert Q^{n+1/2} \Vert_{\infty}\right)T\,, \quad 0\leq l\leq N\,,
\end{equation}
where $T=N\Delta t$. The inequality \eqref{eq_estth} is obtained by applying \eqref{eq_estphi} and Lemma \ref{le_estst} again, together with the following inequality 
\[
\left\Vert \Psi^l\right\Vert_2 \leq \left\Vert \varphi^l\right\Vert_2 + \left\Vert \theta^l\right\Vert_2 \leq \left\Vert \varphi^l\right\Vert_2 + 2\left\Vert \Phi\right\Vert_2\,. 
\]
\end{proof}

\subsubsection{Simplified discrete transparent boundary conditions}

As proposed in \cite{EhAr}, the computational time required by our transparent boundary conditions can be reduced by using the decay of the coefficients $s_j^n$ appearing in \eqref{dtbcNHg}\eqref{dtbcNHd}. This improvement is obtained by replacing the coefficients $s_j^n$ by $0$ after an index $M\geq 1$ chosen big enough to preserve the accuracy of the boundary conditions. In that case, the numerical cost of the convolutions in \eqref{dtbcNHg}\eqref{dtbcNHd} is reduced and the computations are performed as follows: for $1\leq l \leq M+1$ the boundary conditions are given by \eqref{dtbcNHg}\eqref{dtbcNHd} and for $l\geq M+2$, they are given by: 
\begin{multline}\label{dtbcNHg_cut}
\displaystyle \Psi^l_1 - s^0_0\Psi^l_0 = \sum_{k=l-M}^{l-1}s^{l-k}_0\left(\Psi^k_0- \Phi_0e^{-i\omega_0t^k}\right)-\left(\Psi^{l-1}_1-\Phi_1e^{-i\omega_0t^{l-1}}\right) + \Phi_1e^{-i\omega_0t^l} - s^0_0\Phi_0e^{-i\omega_0t^l}
\end{multline}
and
\begin{multline}\label{dtbcNHd_cut}
\displaystyle \Psi^l_{J-1} - s^0_J\Psi^l_J = \sum_{k=l-M}^{l-1}s^{l-k}_J\left(\Psi^k_J-\Phi_Je^{-i\omega_Lt^k}\right) - \left(\Psi^{l-1}_{J-1}-\Phi_{J-1}e^{-i\omega_Lt^{l-1}}\right) + \Phi_{J-1}e^{-i\omega_Lt^l} - s^{0}_J\Phi_Je^{-i\omega_Lt^l}\,.
\end{multline}

\end{document}